\documentclass[AMA,STIX1COL]{WileyNJD-v2}

\articletype{Article Type}%

\received{x}
\revised{x}
\accepted{x}

\usepackage{amsmath}
\usepackage{siunitx}
\usepackage{subfig}
\usepackage{longtable,tabularx,multirow}
\usepackage{cases}

\renewcommand{\vec}[1]{\mbox{\boldmath$#1$}}
\newcommand{\tensor}[1]{\mbox{\boldmath{\ensuremath{#1}}}}

\usepackage{hyperref}
\usepackage[noabbrev]{cleveref}

\begin{document}

\title{On the numerical accuracy in finite-volume methods to accurately capture turbulence in compressible flows}

\author[1]{Emmanuel Motheau*}

\author[2]{John Wakefield}

\authormark{Motheau and Wakefield}

\address[1]{\orgdiv{Center for Computational Sciences and Engineering, Computational Research Division}, \orgname{Lawrence Berkeley National Laboratory}, \orgaddress{1 Cyclotron Rd, MS 50A-3111, Berkeley, CA 94720, \country{USA}}}

\address[2]{\orgdiv{Department of Mathematics}, \orgname{University of Michigan}, \orgaddress{East Hall, 530 Church
Street, Ann Arbor, MI 48109-1043, \country{USA}}}

\corres{*Corresponding author. \email{emotheau@lbl.gov}}

\presentaddress{Lawrence Berkeley National Laboratory, \orgaddress{1 Cyclotron Rd, MS 50A-3111, Berkeley, CA 94720, \country{USA}}}

\abstract[Summary]{
The goal of the present paper is to understand the impact of numerical schemes for the reconstruction of data at cell faces in finite-volume methods, and to assess their interaction with the quadrature rule used to compute the average over the cell volume. Here, third-, fifth- and seventh-order WENO-Z schemes are investigated.  On a problem with a smooth solution, the theoretical order of convergence rate for each method is retrieved, and changing the order of the reconstruction at cell faces does not impact the results, whereas for a shock-driven problem all the methods collapse to first-order. Study of the decay of compressible homogeneous isotropic turbulence reveals that using a high-order quadrature rule to compute the average over a finite-volume cell does not improve the spectral accuracy and that all methods present a second-order convergence rate. However the choice of the numerical method to reconstruct data at cell faces is found to be critical to correctly capture turbulent spectra. In the context of simulations with finite-volume methods of practical flows encountered in engineering applications, it becomes apparent that an efficient strategy is to perform the average integration with a low-order quadrature rule on a fine mesh resolution, whereas high-order schemes should be used to reconstruct data at cell faces.
}



\maketitle

\section{Introduction}

There are many misconceptions in Computational Fluid Dynamics (CFD) pertaining to the accuracy of numerical methods on finite meshes. Indeed, it is often assumed that the theoretical convergence rate of the numerical error of a scheme fully describes the accuracy for a finite mesh (see the discussion in \cite{Wang:2013}). Therefore, one would assume that high-order methods reach a desired solution faster than low-order methods because they would require fewer grid points.

In a previous work, Motheau and Wakefield  \cite{Motheau:2020} investigated this subject in the context of finite-volume methods by performing a fair (i.e.\ with consideration of computational cost) comparison between low (second-order) and high-order (fourth-order or greater) schemes for a set of finite meshes up to what is required to resolve physical scales. It is emphasized that there are a plethora of different numerical strategies to solve partial differential equations, but in the context of CFD, especially for complicated applications encountered in engineering, finite-volume methods are often preferred because they are intrinsically conservative, shock-capturing, and flexible enough to handle both unstructured and structured meshes. Further they are well-suited to use in Adaptive Mesh Refinement (AMR) through re-fluxing across multi-grids to achieve conservation properties.

 The study presented in \cite{Motheau:2020} shows that for a problem involving a smooth solution, the high-order method performs better than the second-order ones and that the theoretical order of convergence rate is retrieved. However, when the solution contains a shock all the methods collapse to a first-order convergence rate. In the context of the decay of compressible homogeneous isotropic turbulence with shocklets, which is representative of the turbulent activity in a realistic compressible flow, the actual overall convergence rate of the methods reduces to second-order. Furthermore, one of the main outcome of this work was to demonstrate that in terms of turbulent spectra, all the numerical methods provide similar results and that virtually the same physical solution can be obtained much faster by refining a simulation with the second-order method rather than running a coarse high-order simulation. 

It becomes clear that evaluating the accuracy of a numerical method for turbulent flows should also consider its spectral dissipation and dispersion properties on meshes of interest in addition to its theoretical formal order of convergence. Moreover, as mentioned in Motheau and Wakefield \cite{Motheau:2020}, the theoretical order of convergence is often not proven for these methods as it is the result of the combination of different numerical schemes together to perform different tasks. In the finite-volume context, numerical integration is performed with two major steps. First, the physical state is reconstructed cell faces in order to evaluate fluxes. Second, a quadrature rule is applied to compute the average at the cell center. The numerical properties of the scheme employed for the reconstruction of data at cell faces may interact with the one employed for the quadrature rule used to compute the average at cell center in complex ways. 

In the work of Motheau and Wakefield  \cite{Motheau:2020}, two different numerical strategies were tested and compared: a second-order Godunov method with PPM~interpolation, as well as the fourth-order finite-volume WENO method proposed by \cite{Titarev:2004}. The Godunov method used to compute fluxes is indeed second-order, but in this work the reconstruction of data at cell faces is  performed either with fourth-order interpolation (as the original PPM method \cite{Colella:2008}) or by fifth-order WENO interpolation \cite{WENO_JS}. In the strategy of \cite{Titarev:2004}, the same fifth-order WENO scheme is employed to evaluate fluxes at faces, but a fourth-order quadrature rule is employed to compute the integral over a cell.

The goal of the present paper is to understand the impact of numerical schemes for the reconstruction of data at cell faces and to assess their interaction with the volume integration method. The present study takes advantage of the hybrid PPM/WENO strategy that has been developed in  \cite{Motheau:2020}. Whereas the whole integration scheme is a second-order Godunov method with PPM interpolation, the reconstruction of data at cell faces can be performed by any WENO scheme. Here, third-, fifth- and seventh-order WENO-Z schemes are investigated. Note that instead of the original WENO method \cite{WENO_JS}, the WENO-Z variant \cite{Borges:2008} is employed because its superior performance and robustness in the context of compressible turbulent flows was demonstrated in \cite{Motheau:2020}. Furthermore, in order to assess the influence of the volume integration rule, results are compared to the fourth-order finite-volume method \cite{Titarev:2004} with fifth-order WENO-Z schemes for reconstruction at cell faces.

Three test cases of increasing complexity are investigated in the present paper. First, the convection of a smooth vortex is considered, followed by the simulation of a classical shock-driven Shu-Osher problem. Finally, the decay of compressible homogeneous isotropic turbulence (HIT) with shocklets is investigated. Results are consistent with the work of Motheau and Wakefield  \cite{Motheau:2020}. On a problem with a smooth solution, the theoretical order of convergence rate for each method is obtained and changing the accuracy of the reconstruction of data at the cell face does not impact the results. Again, for the Shu-Osher problem all the methods collapse to first-order. The HIT problem shows interesting results; whereas all methods present a second-order convergence rate for the physical data, analysis of the turbulent spectra reveals that the choice of a numerical scheme for the reconstruction of data at cells has a significant strong impact on the spectra, and that it does not depend on the volume integration quadrature rule.

Overall, the present study reveals that using a high-order method for the integration rule over a finite-volume cell does not improve the spectral accuracy. The choice of the numerical method to reconstruct data at cell faces is however critical to correctly capture turbulent spectra. Because high-order methods for the integration over a cell require a significant number of evaluations of face data reconstruction, it becomes clear that an efficient strategy for finite-volume methods is to choose a low-order method to perform the average integration over the cells on a fine mesh discretization, and that high-order accurate schemes should be used to reconstruct data at cells face.

\section{Numerical methods}

The present study is performed with the \textbf{PeleC} code, which is a second-order AMR finite-volume solver for reacting and non-reacting fluid simulations with complex geometry and support for Lagrangian spray particles. The simulations performed in the present paper only use a fraction of the capability of the software, namely the Godunov-based hybrid PPM/WENO integration procedure on a single level mesh grid. Note that \textbf{PeleC} is part of the \textbf{Pele Suite} of codes, which are publicly available and may be freely downloaded\footnote{\url{https://amrex-combustion.github.io/}}, and that all the test cases investigated in the present paper are available from the \textbf{PeleC} distribution and can be reproduced. Also, the spectral and temporal analysis performed in \cref{subsec:compressible_HIT} can be reproduced with the tools provided in the \textbf{PeleAnalysis} repository of the \textbf{Pele Suite}. Furthermore, in order to provide comparisons with a high-order finite-volume integration scheme, the \textbf{RNS} code is also employed. \textbf{RNS} implements the Adaptive Multi-Level Spectral Deferred Correction (AMLSDC) method, which is fourth-order in time \cite{Emmett:2018}, as well as the fourth-order finite-volume WENO method \cite{Titarev:2004}.

Whereas both codes were initially developed for the simulation of combustion problems and solve the multicomponent reacting Navier-Stokes equations, only non-reacting problems with no specific mixture are investigated in the present study. Consequently, both codes solve the same system of simplified equations, which is given by
\begin{align}
\frac{\partial \rho}{\partial t} + \nabla \cdot (\rho
    \vec{u})= { } & 0, \label{eq:NS:mass} \\
\frac{\partial \rho \vec{u}}{\partial t} + \nabla \cdot (\rho
    \vec{u} \otimes \vec{u}) + \nabla p= { } & \nabla \cdot \tensor{\tau}, \label{eqs:NS:momentum} \\
\frac{\partial \rho E}{\partial t} + \nabla \cdot [(\rho E + p)
  \vec{u}] = { } & \nabla \cdot (\lambda \nabla T) + \nabla \cdot
  \left(\tensor{\tau} \cdot \vec{u}\right), \label{eqs:NS:energy}
\end{align}
where $\rho$ is the density, $\vec{u}$ is the velocity, $p$ is the pressure, $E = e + \mathbf{u} \cdot \mathbf{u} / 2$ is the total energy, $T$ is the temperature and $\lambda$ is the thermal conductivity. The viscous stress tensor is given by
\begin{equation}
\tensor{\tau} = \eta (\nabla \vec{u} + (\nabla \vec{u})^T) + (\varsigma - \frac{2}{3} \eta ) (\nabla \cdot \vec{u}) \mathbf{I},
\end{equation}
where $\eta$ and $\varsigma$ are the shear and bulk viscosities.

The system is closed by an equation of state (EOS) that specifies $p$ as a function of $\rho$ and $T$.  An ideal gas mixture for the EOS is assumed:
\begin{equation}
 p = \rho T \mathfrak{R},
 \label{eqn:eos}
\end{equation} 
where $\mathfrak{R}$ is the specific gas constant. Here we set $C_p$ and $C_v$ the heat capacity at constant pressure and volume, respectively, to follow an ideal gas law proportional to the ratio of the specific heats $\gamma$ so that \cref{eqn:eos} is equivalent to the following relation:
\begin{equation}
e = \frac{p}{\gamma - 1} \, \rho 
\end{equation}
where $e$ is the specific internal energy and $\gamma$ is set to $\gamma=1.4$.

To keep the present paper as simple and concise as possible, the spatial and temporal integration schemes are not presented here, but all the details can be found in the reference \cite{Motheau:2020} for the hybrid PPM/WENO method implemented in the \textbf{PeleC} code, as well as \cite{Titarev:2004} and \cite{Emmett:2018} for the fourth-order finite-volume WENO method and the fourth-order AMLSDC implemented in the \textbf{RNS} code, respectively.

The present paper investigates the reconstruction of data at cell faces with WENO-Z schemes with three different orders of accuracy. Recall that for a given cell $i$, the principle of a WENO method is to provide a high-order approximation of the variable $q$ interpolated on the left and right side of a face, denoted $\hat{q}^L_{i+\frac{1}{2}}$ and $\hat{q}^R_{i-\frac{1}{2}}$. Thus, a $(2r-1)$-th order polynomial approximation of $\hat{q}_{i\pm\frac{1}{2}}$ is constructed as a convex combination of interpolants yielded by varied stencils $ k $  of the values denoted $\hat{q}^k_{i\pm\frac{1}{2}}$:
\begin{equation}
\hat{q}_{i\pm\frac{1}{2}} = \sum_{k=0}^{r-1} \omega_k \hat{q}^k_{i\pm\frac{1}{2}} \label{eqn:WENO_approx}
\end{equation}
where $ \omega_k $ are chosen with the goal of selecting the smoothest such combination.  Consequently, a third-order WENO scheme will lead to $k=0,1,2$, while a fifth- and a seventh-order scheme will lead to $k=0,1,2,3$ and $k=0,1,2,3,4$, respectively.

In \cref{eqn:WENO_approx}, $\omega_k$ are the non-linear weights balancing the contribution of each stencil, and the challenge is to find the best values to capture shocks the most accurately while preserving the resolution of the spectrum of a solution. The weights $\omega_k$ are defined as
\begin{equation}
\omega_k = \frac{\alpha_k}{\sum^{r-1}_{l=0}\alpha_l},\hspace{1cm} \alpha_k = \frac{d_k}{\left( \beta_k + \epsilon \right)^a},
\label{eqn:weights_WENO}
\end{equation}
where $d_k$ are the so-called optimal weights because they reconstruct the upstream $(2r-1)$-th order central scheme, $\beta_k$ are the smoothness indicators, $\alpha_k$ are the un-normalized weights, and $\epsilon$ is a parameter set to avoid a division by zero. The parameter $a$ controls the adaptation rate. In order to keep the present paper simple and the most concise as possible, the numerical expressions of $d_k$, $\beta_k$ and $\hat{q}^k_{i\pm\frac{1}{2}}$ will not be reproduced here, but can be found in the following references: \cite{WENO_JS} for the third- and fifth-order schemes, and in \cite{BALSARA:2000} for the seventh-order scheme. According to \cite{Arshed:2013}, a large value of $a$ leads to unnecessary dissipation in smooth regions of the flow. In the present study, the parameter is set to $a=2$ for all the test cases and $\epsilon$ is set to $\epsilon=10^{-40}$, whatever the order of accuracy of the WENO scheme.

A well-known issue with the original WENO method is that the smoothness indicators $\beta_k$ employed to compute the weights $\omega_k$ fail to recover the maximum order of the scheme at critical points when the derivatives of flux function vanish. Borges \textit{et al.} \cite{Borges:2008} propose a different approach to overcome the issues of the original WENO method by acting directly on the smoothness indicator $\beta_k$ with a very simple formulation. The so-called WENO-Z method is given by
\begin{equation}
\omega_k^{(\rm z)} = \frac{\alpha_k^{(\rm z)}}{\sum_{i=0}^2\alpha_i^{(\rm z)}}, \hspace{.5cm} \text{with} \hspace{.5cm} \alpha_k^{(\rm z)} = d_k\left(1 + \frac{\tau}{\beta_k + \epsilon} \right)^a,
\end{equation}
where $\tau$ is the absolute value of the difference between $\beta_0$ and $\beta_{r-1}$, such that $\tau = |\beta_0 - \beta_1|$  for the third-order WENO-Z scheme, and $\tau = |\beta_0 - \beta_2|$ and $\tau = |\beta_0 - \beta_3|$ for the fifth- and seventh-order schemes, respectively.

\section{Results}

In this section, different order of accuracy of the WENO-Z reconstruction scheme at faces are implemented within the second-order hybrid PPM/WENO method, and compared to the fourth-order finite-volume WENO method with fifth-order WENO-Z reconstruction at faces edges. \Cref{tab:methods_summary} summarizes the different combination of numerical schemes investigated in the present study.

\begin{center}
\begin{table}[ht]%
\centering
\caption{Summary of the combination of numerical methods investigated for volume integration and data reconstruction at faces.}\label{tab:methods_summary}%
\begin{tabular*}{280pt}{@{\extracolsep\fill}lcc@{\extracolsep\fill}}
\toprule
\textbf{Software} & \textbf{Finite volume integration}  & \textbf{Face reconstruction}  \\
\midrule
\textbf{RNS} & $4$th-order WENO  & $5$th-order WENO-Z   \\
\midrule
\textbf{PeleC} & $2$nd-order PPM  & $3$rd-order WENO-Z \\ 
 &  &  $5$th-order WENO-Z \\
 &  & $7$th-order WENO-Z   \\
\bottomrule
\end{tabular*}
\end{table}
\end{center}


Three different test cases are considered in this work:
\begin{itemize}
	\item the convection of a smooth compressible vortex in a periodic domain is a simple, standard, and efficient test case to recover and highlight the theoretical order of accuracy of a numerical method;
	\item the classic Shu-Osher problem, which represents the extreme opposite of the smooth vortex test case. The Shu-Osher problem is very difficult to solve numerically, because a shock wave is propagating in an oscillating entropy field, and the challenge is to capture the shock while resolving the phase and amplitude of the fluctuating entropy. As it will be shown, all the methods perform correctly, but for all of them the rate of convergence collapses to first-order;
	\item the decay of compressible homogeneous isotropic turbulence in the presence of eddy shocklets. This test case can be viewed as a combination of the two previous test case, because it contains both shocks and discontinuities, as well as smooth turbulence structures that lie in a large-bandwidth turbulent spectrum. More specifically, this test case is representative of flows that are encountered in practical CFD applications (see  \cite{Motheau:2014b} for an example).
\end{itemize}

\subsection{2D convection of a smooth compressible vortex}
\label{subsec:COVO}

The following test case consists of the convection of a 2D compressible vortex in a periodic domain so as to accumulate numerical errors from the discretization schemes. This test case exhibits a smooth solution and weak compressibility effects, which allows numerical schemes to perform near their asymptotic limit \cite{Motheau:2017aa,Motheau:2020}. Simulations are performed with increasing mesh resolution and the time-step is computed based on the mesh resolution via a constraint on the CFL number, set to $0.7$. At the end of a simulation, convergence is measured using the $\mathcal{L}^1$-norm of the difference of the $x$-velocity between the final computed solution and the initial solution, and is denoted $\varepsilon_{u}$.
%
%

The configuration is a single vortex superimposed on a uniform diagonal flow field along the $x$- and $y$-directions. The stream function $\Psi$ of the initial vortex is given by
\begin{equation}
\Psi = \Gamma \exp\left(\frac{-r^2}{2 R_v^2} \right),
\end{equation}
where $r=\sqrt{\left(x-x_v\right)^2 + \left(y-y_v\right)^2}$ is the radial distance from the center of the vortex located at $\left[x_v, y_v\right]$ (the center of the domain), while $\Gamma$ and $R_v$ are the vortex strength and radius, respectively. The velocity field is then defined as
\begin{equation}
u = u_{0} + \frac{\partial  \Psi}{\partial y}, \hspace{1cm} v = u_{0} -\frac{\partial  \Psi}{\partial x}.
\end{equation}
The initial pressure field is expressed as
\begin{equation}
p\left(r\right) = p_{\rm ref} \exp\left(-\frac{\gamma}{2} \left(\frac{\Gamma}{c R_v} \right)^2 \exp \left( -\frac{r^2}{R_v^2} \right) \right),
\end{equation}
and the corresponding density field is given by
\begin{equation}
\rho\left(r\right) = \frac{p\left(r\right)}{\mathcal{R} T_{\rm ref}}, 
\end{equation}
where $T_{\rm ref}$ is assumed constant and $\gamma=1.4$.

The computational domain is a square of dimension $L=0.01$~m. The reference temperature $T_{\rm ref}$ and pressure $p_{\rm ref}$ are set to $300$~K and $101320$~Pa, respectively. Parameters are set to $\Gamma = 0.11$~m$^2$/s and $R_v = 0.1 L$. The initial flow velocity is $u_{0} = 100$~m/s. In the present test case, only the Euler equations  are solved; there is no diffusion.

The simulations are performed over a physical time of $5$~ms, corresponding to $5$ flow through times (FTT), in order to accumulate numerical errors from the spatial discretization.

\begin{figure}[th!]
\centering
\includegraphics[width=0.5\textwidth]{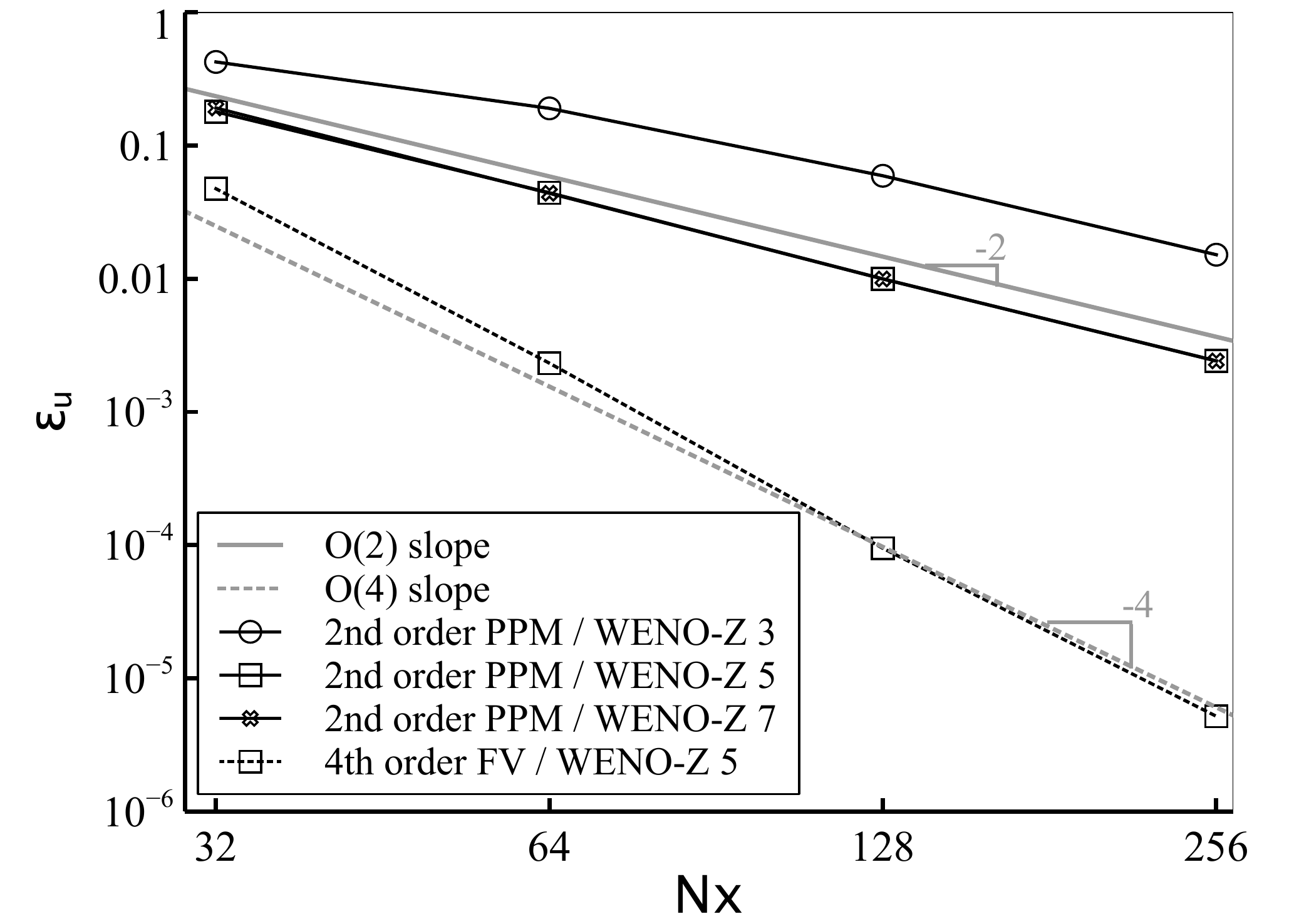}
\caption{Convection of a vortex, evolution of the $\mathcal{L}^1$ norm of the error of the $x-$velocity for different mesh size $N_x$.}
\label{fig:COVO_error}
\end{figure}

Results are shown in \cref{fig:COVO_error}. The solid and dotted grey lines represent  second- and fourth-order slopes, respectively. As expected, because the solution is smooth, all the numerical methods show convergence at their maximum theoretical order of accuracy. The fourth-order finite-volume method is indeed converging towards a fourth-order asymptotic limit, whereas the hybrid PPM/WENO method is second-order. The present simulations show that the order of accuracy of the reconstruction of data at face edges does not impact the overall accuracy of the method. Because the solution is smooth, numerical errors at edges are very small and the errors are caused by the volume integration performed by the second-order PPM method.

\subsection{Shock-driven test case: the Shu-Osher problem}
\label{subsec:Shu_Osher}

The so-called Shu-Osher test case simulates the one-dimensional propagation of a normal shock wave interacting with a fluctuating entropy wave, generating a flow field containing both small scale structures as well as discontinuities. The initial conditions are given by:
\begin{equation}
    \left( \rho, u, p \right) = \begin{cases}
        \left(3.857143, 2.629369, 10.3333 \right), & \text{if}\;x \leqslant 1,  \\
        \left( 1+ 0.2 \sin \left(5x\right),0,1\right), & \text{otherwise}.
    \end{cases}
\end{equation}
The length of the computational domain is $x \in [0,10]$ and the solution is advanced in time to $t=1.2$. For all numerical methods investigated, the mesh is progressively refined from $N_x=256$ to $N_x = 2048$. The convergence is measured using the $\mathcal{L}^1$-norm of the difference in density between the final computed solution and a reference solution defined to be the solution computed with the fourth-order finite-volume WENO method and with a very fine mesh $N_x =32768$. In all simulations the CFL number is set to $0.5$.

\begin{figure}[thp]
\centering
\subfloat[Full domain]{\label{fig:Shu_Osher_WENO_PPM_256_a}\includegraphics[width=0.4\textwidth]{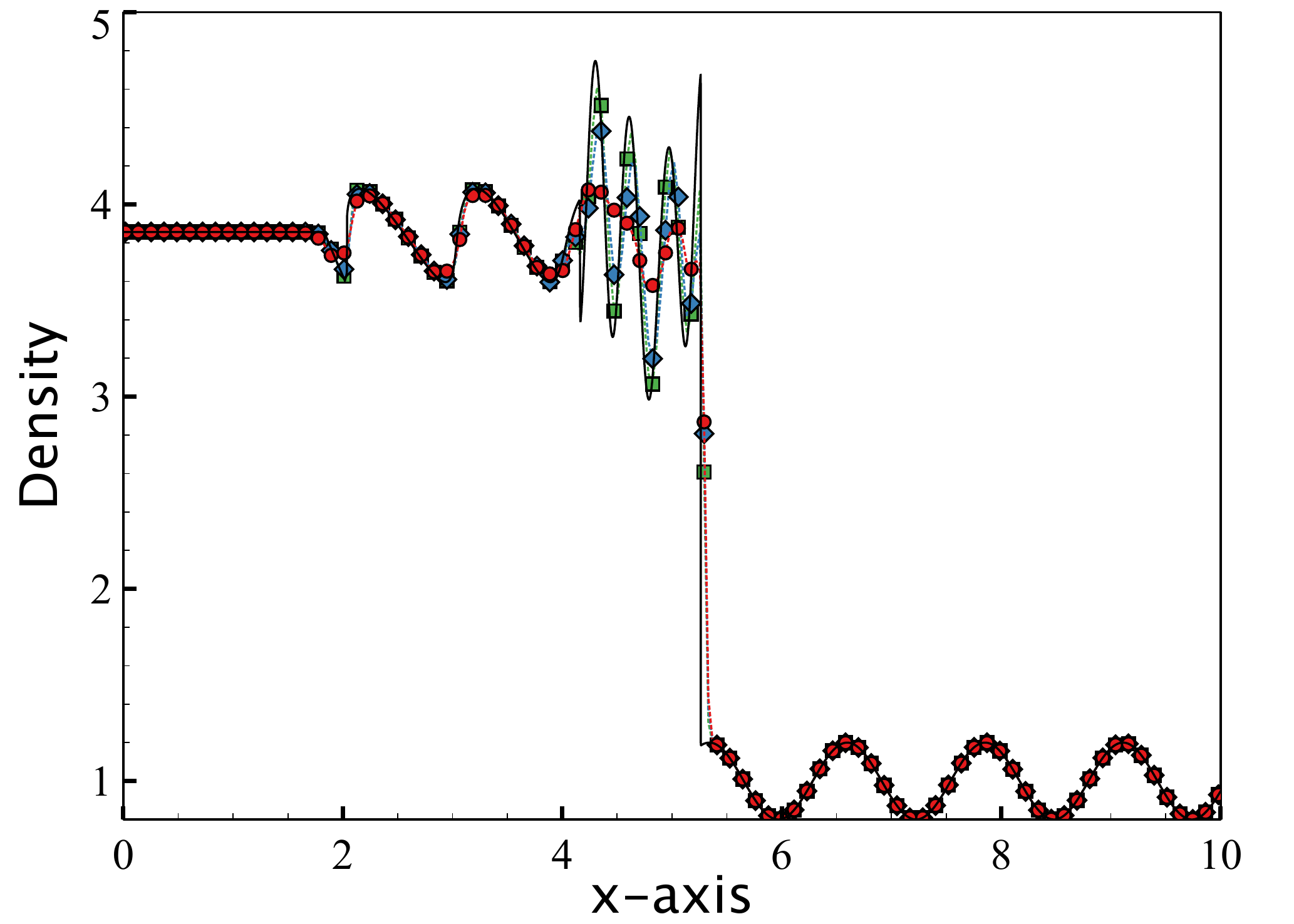}}
\subfloat[Zoom]{\label{fig:Shu_Osher_WENO_PPM_256_b}\includegraphics[width=0.4\textwidth]{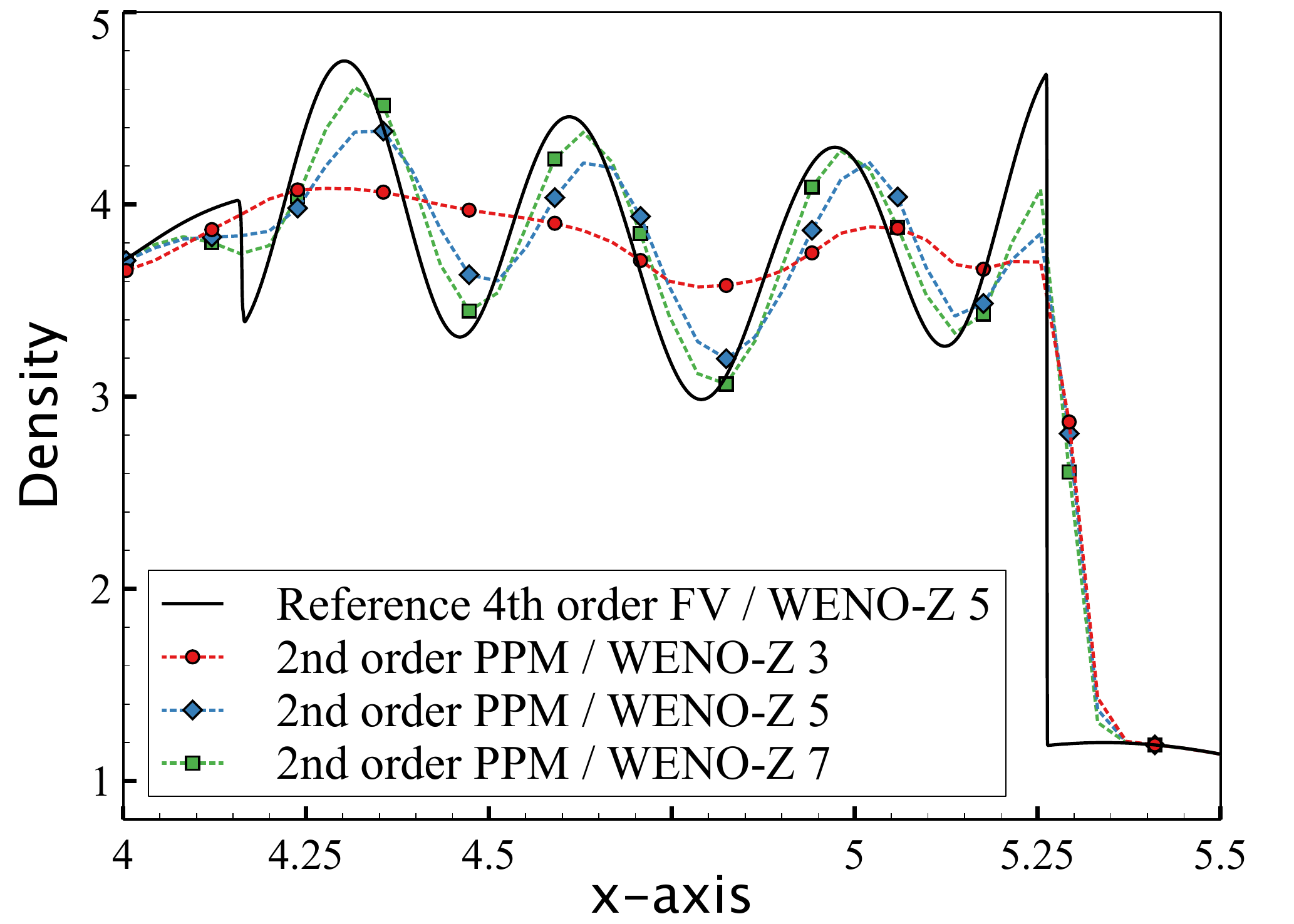}}
\caption{Shu-Osher test case: profile of density for $N_x=256$.}
\label{fig:Shu_Osher_WENO_PPM_256}
\end{figure}

\begin{figure}[thp]
\centering
\subfloat[Full domain]{\label{fig:Shu_Osher_WENO_PPM_512_a}\includegraphics[width=0.4\textwidth]{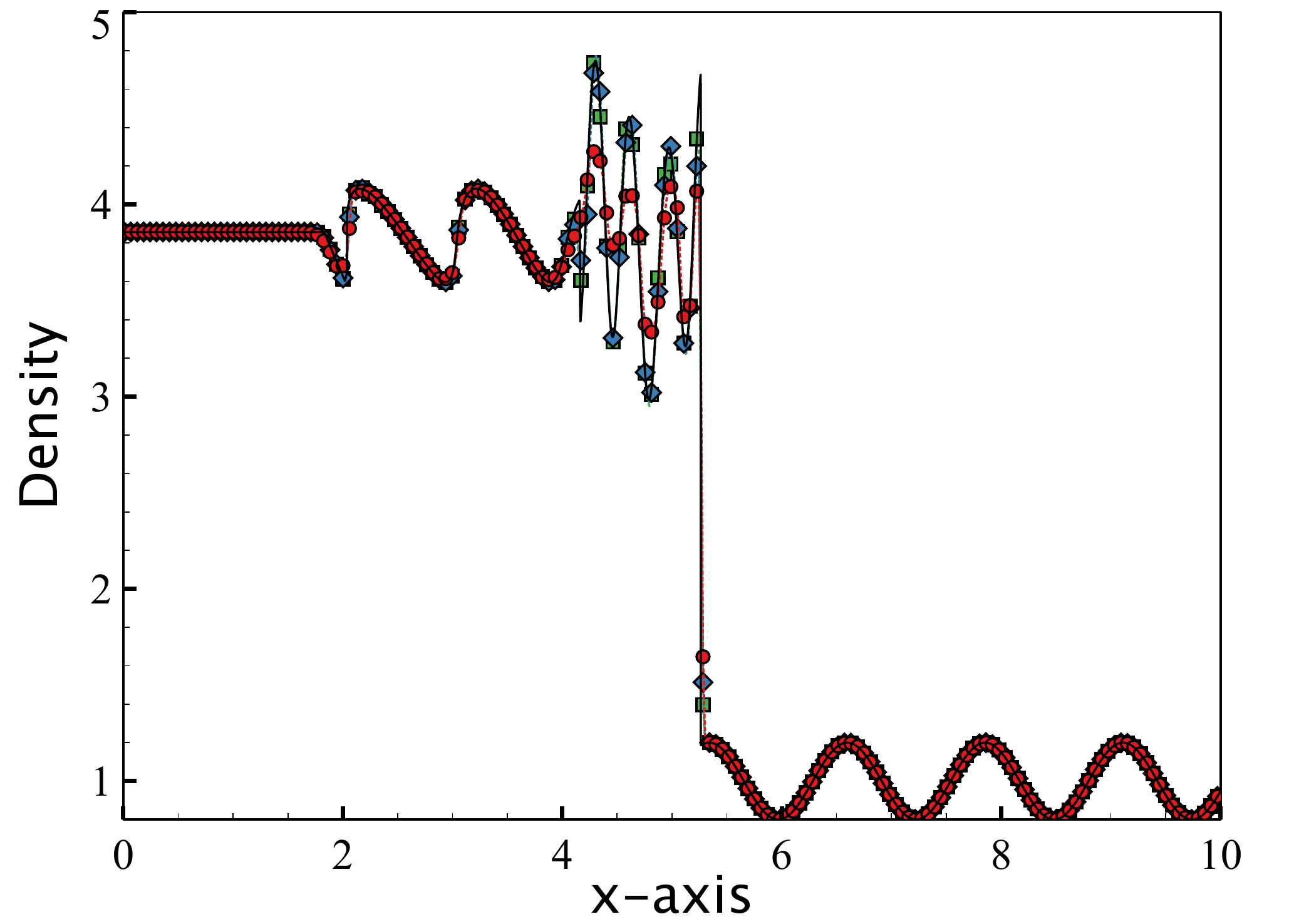}}
\subfloat[Zoom]{\label{fig:Shu_Osher_WENO_PPM_512_b}\includegraphics[width=0.4\textwidth]{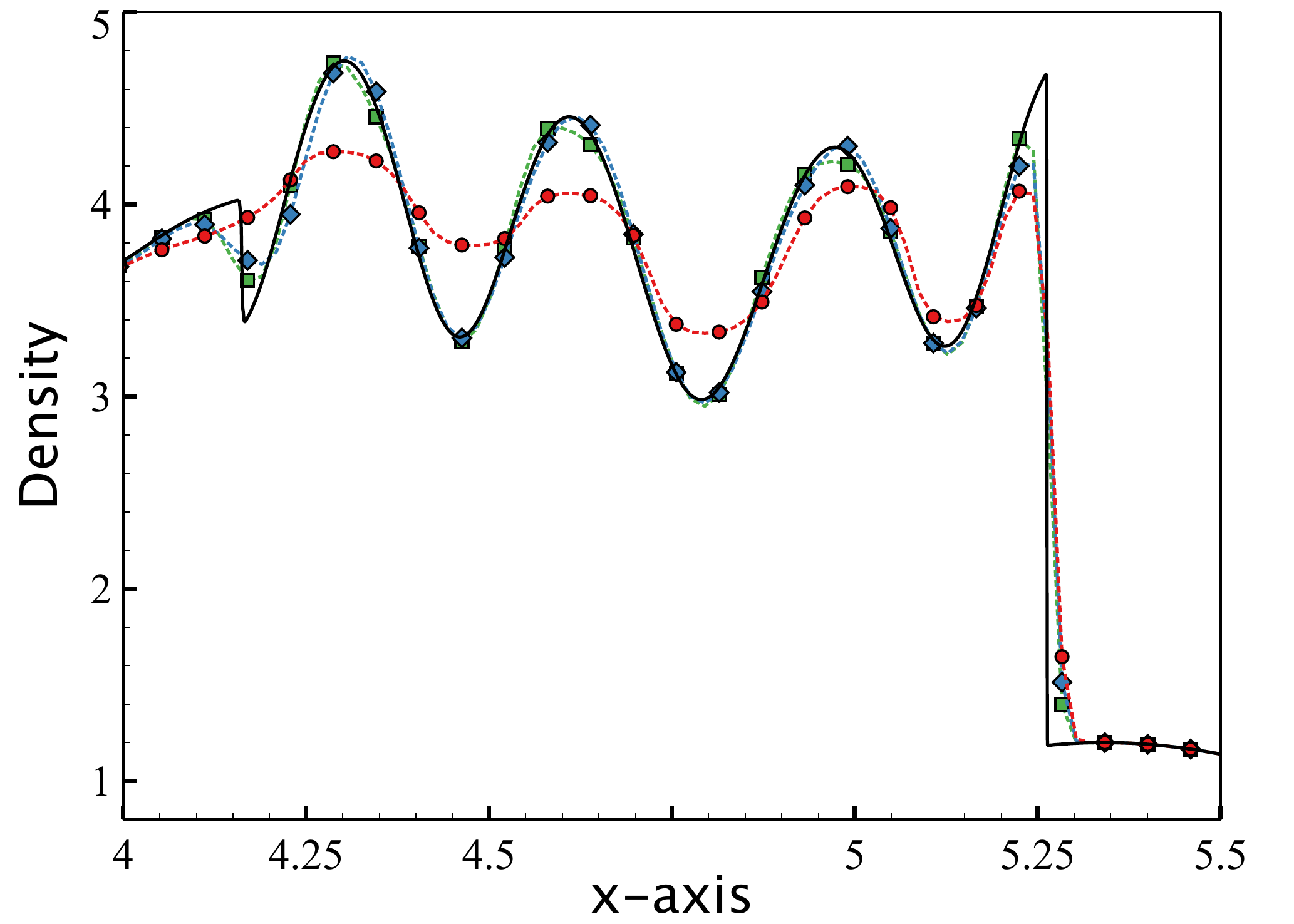}}
\caption{Shu-Osher test case: profile of density for $N_x=512$.}
\label{fig:Shu_Osher_WENO_PPM_512}
\end{figure}

\begin{figure}[thp]
\centering
\subfloat[Full domain]{\label{fig:Shu_Osher_WENO_PPM_1024_a}\includegraphics[width=0.4\textwidth]{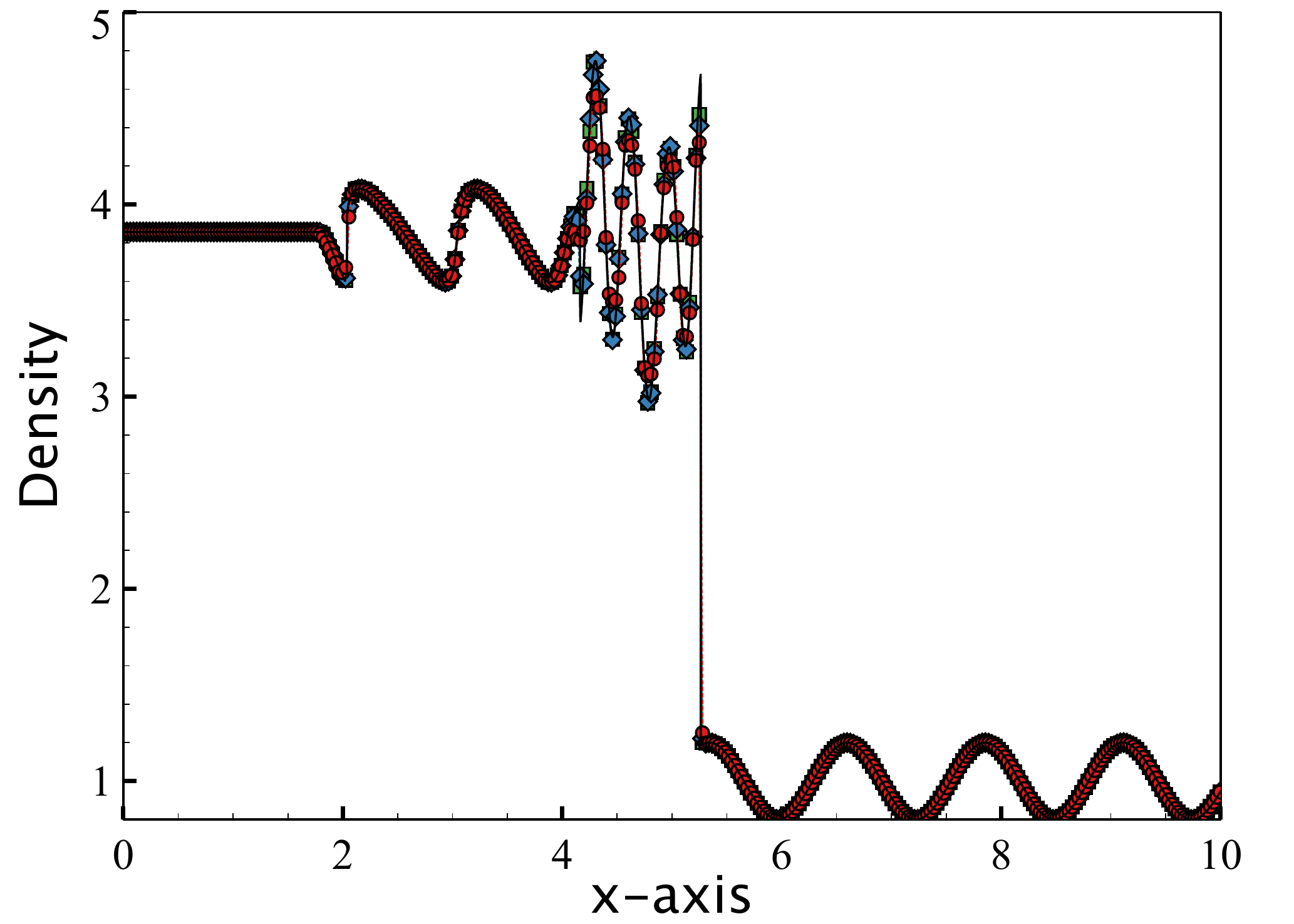}}
\subfloat[Zoom]{\label{fig:Shu_Osher_WENO_PPM_1024_b}\includegraphics[width=0.4\textwidth]{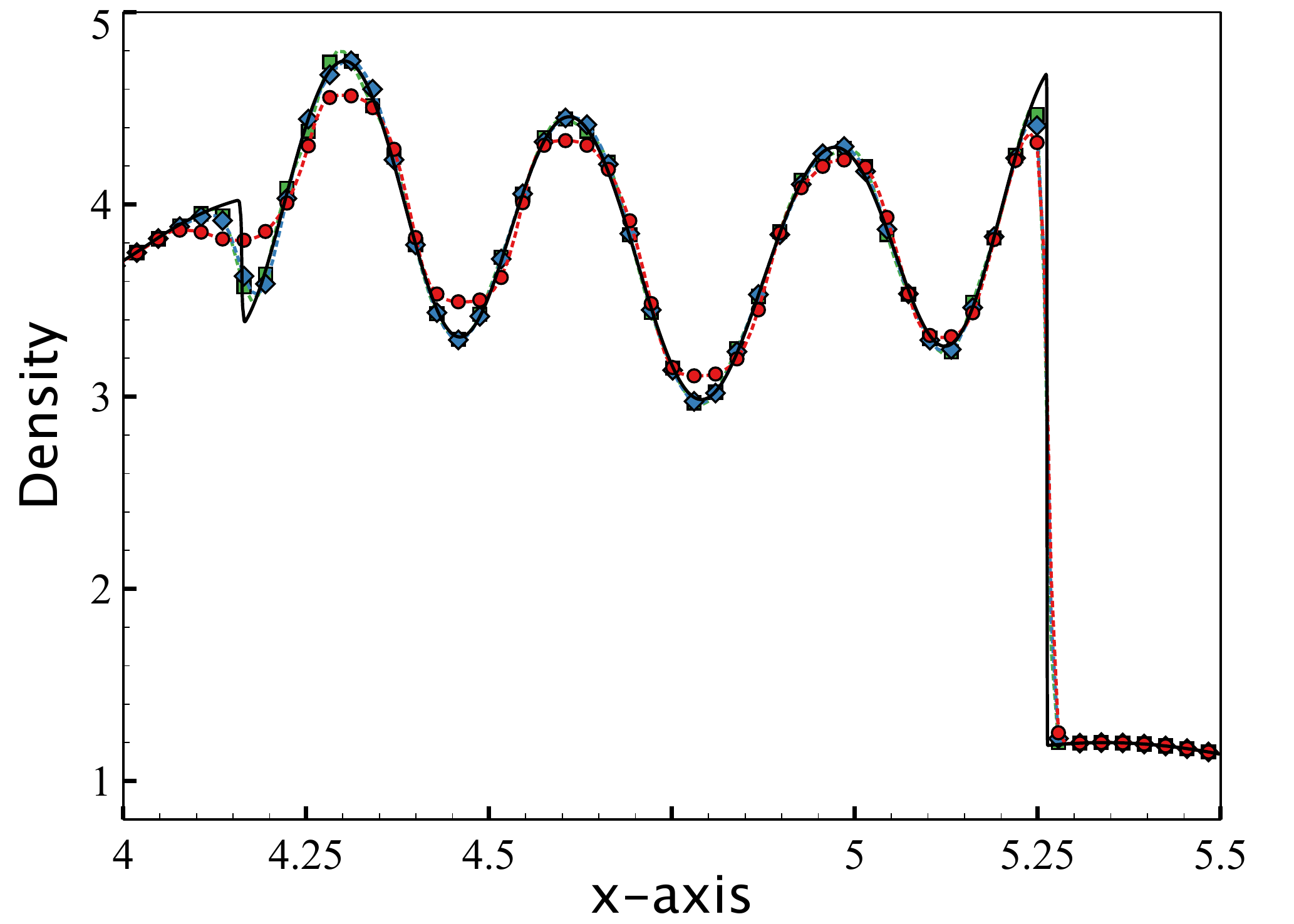}}
\caption{Shu-Osher test case: profile of density for $N_x=1024$.}
\label{fig:Shu_Osher_WENO_PPM_1024}
\end{figure}

\begin{figure}[thp]
\centering
\includegraphics[width=0.5\textwidth]{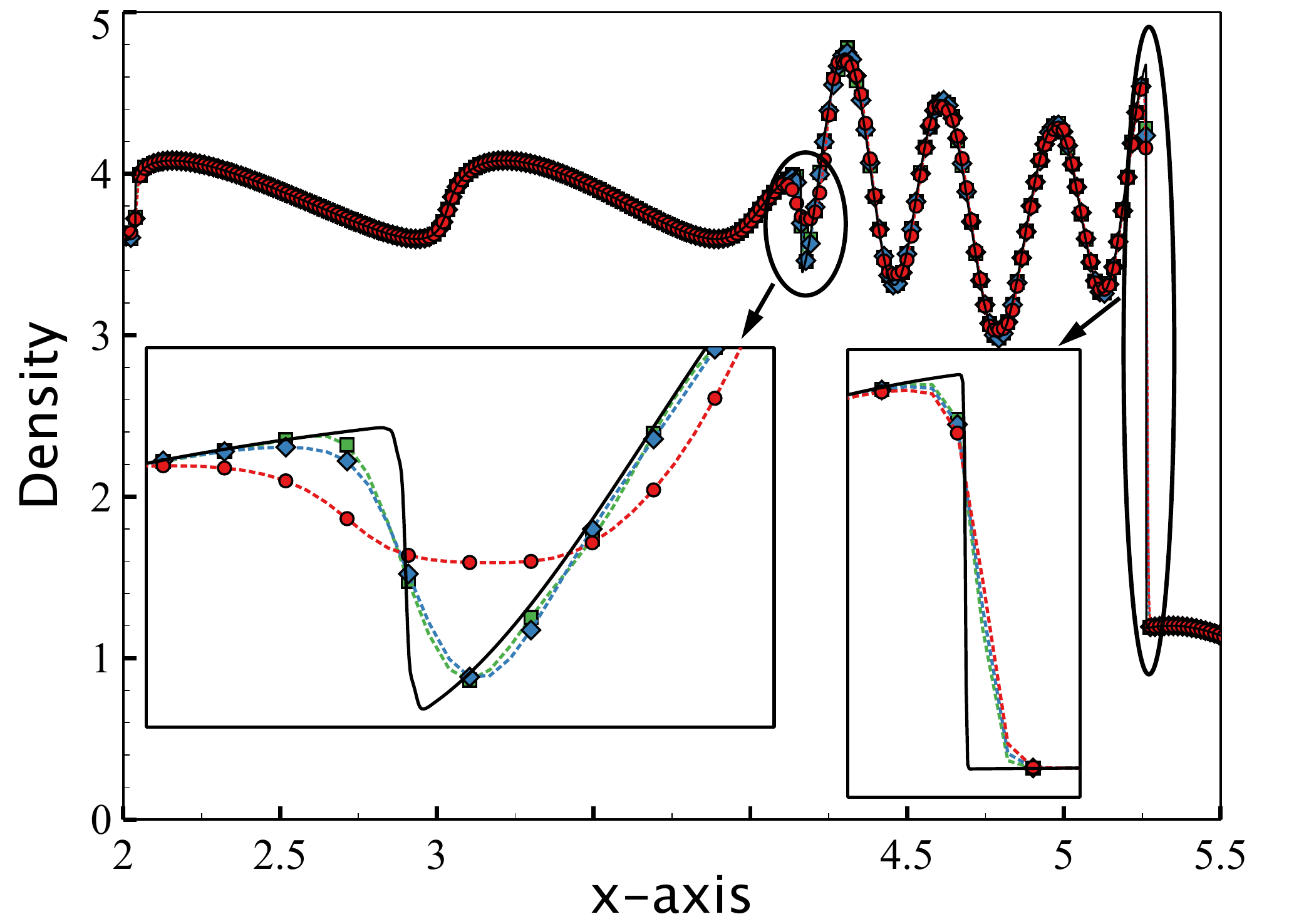}
\caption{Shu-Osher test case: profile of density for $N_x=2048$.}
\label{fig:Shu_Osher_WENO_PPM_2048}
\end{figure}

The density field at $t=1.2$ computed with $N_x = 256, 512, 1024$, and $2048$ is shown in~\cref{fig:Shu_Osher_WENO_PPM_256,fig:Shu_Osher_WENO_PPM_512,fig:Shu_Osher_WENO_PPM_1024,fig:Shu_Osher_WENO_PPM_2048}, respectively. In these figures, the black line is the reference solution, while the red circle, blue diamond and green square symbols represent the solution computed with the second-order hybrid PPM/WENO method and with the reconstruction at face performed with the third-, fifth- and seventh-order WENO-Z schemes, respectively (see legend in~\cref{fig:Shu_Osher_WENO_PPM_256_b}). Note that the panels (a) and (b) in~\cref{fig:Shu_Osher_WENO_PPM_256},~\cref{fig:Shu_Osher_WENO_PPM_512}~and~\cref{fig:Shu_Osher_WENO_PPM_1024} present the full domain and an inset detailing selections of the domain, respectively, while \cref{fig:Shu_Osher_WENO_PPM_2048} shows only that selection of the domain. Note there is no relation between the symbols and the number of grid points. Several symbols have been removed from the figures for clarity purposes.

As it can be seen in \cref{fig:Shu_Osher_WENO_PPM_256}, for a coarse resolution of the mesh ($N_x=256$) the accuracy of the reconstruction of the data at edge faces is crucial, whereas the whole order of accuracy of the integration is theoretically second-order. The third-order WENO-Z scheme performs poorly and totally miss the wave of density as well as the shock profile. Both fifth- and seventh-order schemes are able to capture the overall correct shape of the wave, and the latter one provides a better estimation of the amplitude. In the interest of clarity, the solution for this mesh size with the fourth-order finite volume method has been omitted but can be found in our previous work \cite{Motheau:2020}. It is emphasized that such a solution is similar to the one computed with the hybrid PPM/WENO and the fifth-order WENO-Z scheme (blue curve in \cref{fig:Shu_Osher_WENO_PPM_256_b}).  

All methods, regardless of the order of accuracy with which values at cell faces are reconstructed, provide acceptable results for a refinement of the mesh by a factor $2$ ($N_x=512$), as shown in \cref{fig:Shu_Osher_WENO_PPM_512}. The third-order WENO-Z underestimates the amplitudes of the wave and the shock, but the shape is consistent with the reference solution. Both fifth- and seventh-order schemes better capture the oscillations trailing the shock, with an advantage to the seventh-order scheme as expected. As shown in \cref{fig:Shu_Osher_WENO_PPM_1024} and \cref{fig:Shu_Osher_WENO_PPM_2048}, the fifth- and seventh-order WENO-Z schemes provide nearly identical results very close to the reference solution and the third-order WENO-Z scheme underestimates slightly the amplitude of the wave. Interestingly, the third-order  scheme fails to recover the strong gradient located at $x \approx 4.25$, even with such a fine resolution.

\begin{figure}[thp]
\centering
\includegraphics[width=0.5\textwidth]{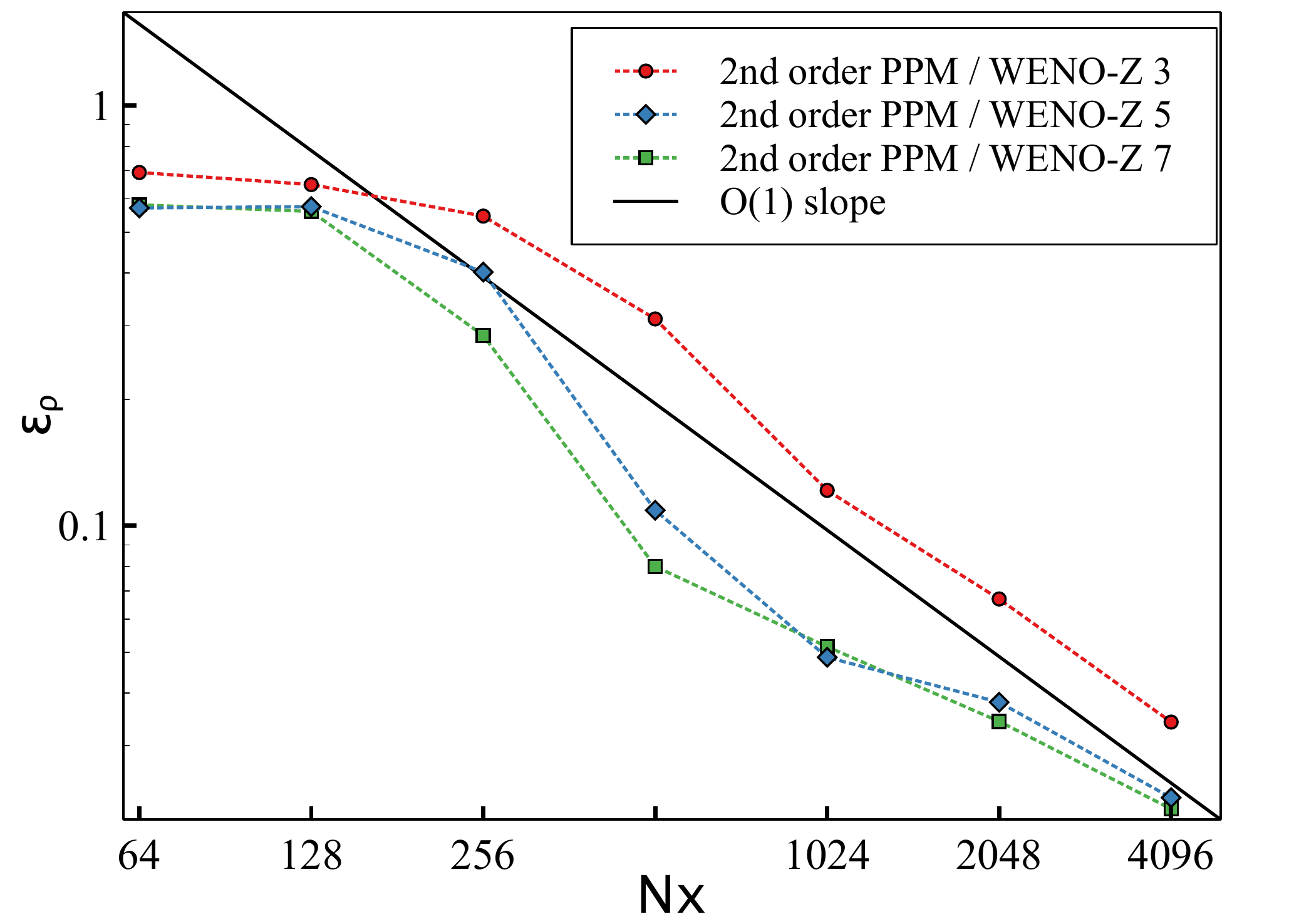}
\caption{Shu-Osher test case: $\mathcal{L}^1$-norm of the error on the density.}
\label{fig:Shu_Osher_convergence_study}
\end{figure}

The evolution of the $\mathcal{L}^1$-norm of the error on the density (denoted $\varepsilon_{\rho}$) is depicted in \cref{fig:Shu_Osher_convergence_study}, and the global convergence rate computed with a best-fitting curve method is reported in \cref{tab:Shu_Osher_convergence_rate}. Whereas rising the order of accuracy of the reconstruction scheme at faces provides less numerical error, which is consistent with the analysis performed above, all the methods investigated collapse to a first-order convergence rate.

Overall, the present study suggests that the theoretical second-order accuracy of the spatial integration method, which was clearly demonstrated for a smooth solution in the previous section \cref{subsec:COVO}, has a significantly smaller impact in the present shock-containing case. It is shown here that the accuracy of the reconstruction of the data at cell faces is crucial for accuracy in smooth regions, though this order of accuracy does not appear in the numerical analysis.

In the following section, a more realistic three-dimensional compressible turbulent flow is simulated to investigate the impact of the numerical accuracy of the different numerical schemes, as well as their effective cost in terms of mesh resolution, when both shocks and small turbulence structures interact in the same domain.

\begin{table}[ht!]%
\centering
\caption{Shu-Osher test case: convergence rate of the $\mathcal{L}^1$-norm of the error on the density.}\label{tab:Shu_Osher_convergence_rate}%
\begin{tabular*}{140pt}{@{\extracolsep\fill}ll@{\extracolsep\fill}}
\toprule
\textbf{WENO-Z scheme} & $\mathcal{O}\left( \epsilon_\rho\right)$ \\
\midrule
$3$rd-order  & $0.90$  \\
$5$th-order  &  $0.95$ \\
$7$th-order  & $0.99$  \\ 
\bottomrule
\end{tabular*}
\end{table}


\subsection{Three-dimensional isotropic compressible turbulence decay}
\label{subsec:compressible_HIT}

The present test case consists on the simulation of the decay of a compressible isotropic turbulent field with the presence of eddy shocklets. Originally a physical study of turbulence in the work of Lee \textit{et al.} \cite{Lee:1991}, this framework of simulations has been employed in \cite{Johnsen:2010,Motheau:2020} to study the properties of numerical schemes to capture turbulence spectra and the decay of physical quantities.

The simulation is controlled by two non-dimensional parameters: the turbulent Mach number
\begin{equation}
M_{t,0} = \frac{\sqrt{<\vec{u}_0 \cdot \vec{u}_0>}}{c_0}
\label{eq:mach_turbulent}
\end{equation}
where $c_0$ is the sound speed in the initial solution, and the Taylor-scale Reynolds number defined as
\begin{equation}
Re_{\psi,0} = \frac{\rho_0  \psi_0 u_{\rm rms,0}}{\eta_0}
\label{eq:reynolds_number}
\end{equation}
where 
\begin{equation}
u_{{\rm rms},0} = \sqrt{\frac{<\vec{u}_0 \cdot \vec{u}_0>}{3}}, \hspace{.5cm} \psi_0 = \frac{2}{k_0}.
\end{equation}
The $\vec{u}_0$ term is a solenoidal velocity field constructed to satisfy the following relation
\begin{equation}
E\left(k\right) \sim k^4 \exp\left(-2\left(k/k_0 \right)^2 \right), \hspace{.5cm} \frac{3 u^2_{{\rm rms},0}}{2} = \frac{<\vec{u}_0 \cdot \vec{u}_0>}{2} = \int_0^{\infty} E\left( k \right){\rm d} k 
\end{equation}
where $k_0$ is the most energetic wavenumber and is set to $k_0=4$ in the present study. 

The study presented here is identical to the one in \cite{Motheau:2020}. In order to allow weak shock waves to develop spontaneously from the turbulent motion, the turbulent Mach and Reynolds numbers are set to $M_{t,0} = 0.6$ and $Re_{\lambda,0}=100$, respectively. This allows numerical convergence for relatively coarse mesh grids so as to keep the computational cost reasonable. Once $M_{t,0}$ and $Re_{\psi,0}$ are set, $u_{{\rm rms},0}$ can be deduced from \cref{eq:mach_turbulent} with the known sound speed, and the viscosity $\eta_0$ can be deduced from  \cref{eq:reynolds_number} and is held constant throughout the simulation. Moreover, a constant thermal conductivity is set according to 
\begin{equation}
\lambda_0 = \frac{\eta_0 C_p}{Pr}
\end{equation}
where $C_p$ is the specific heat capacity, set to $C_p = 1.173$~kJ/kg.K and the Prandtl number $Pr$ is set to $Pr=0.71$. Finally, the initial temperature and pressure in the flow are set to $T_0 = 1200$~K and $p_0=1$~atm.

In the present study, an initial velocity field is generated on a grid of $N_x = 512$ and employed for all the simulations, which are performed over a non-dimensional time set to $t/\tau=4$ where $\tau = \psi_0/u_{{\rm rms},0}$. Several mesh resolutions are investigated: $N_x=64$, $N_x=128$, $N_x=256$ and $N_x=512$, and the CFL number is kept constant at $0.5$. 

In order to assess the impact of the accuracy of numerical schemes for volume integration and face reconstruction, a reference solution is generated with the very high-order code \textbf{SMC} \cite{Emmett:2014}  that employs eighth-order accurate centered finite-difference schemes for the spatial discretization, and a fourth-order Runge-Kutta algorithm for the time advancement. A convergence study for the reference solution is presented in \cite{Motheau:2020}. This reference solution will be depicted with a black solid line in the remainder of the paper.

\begin{figure}[hbt!]
\centering
\subfloat[Kinetic Energy]
{
\includegraphics[width=.45\textwidth]{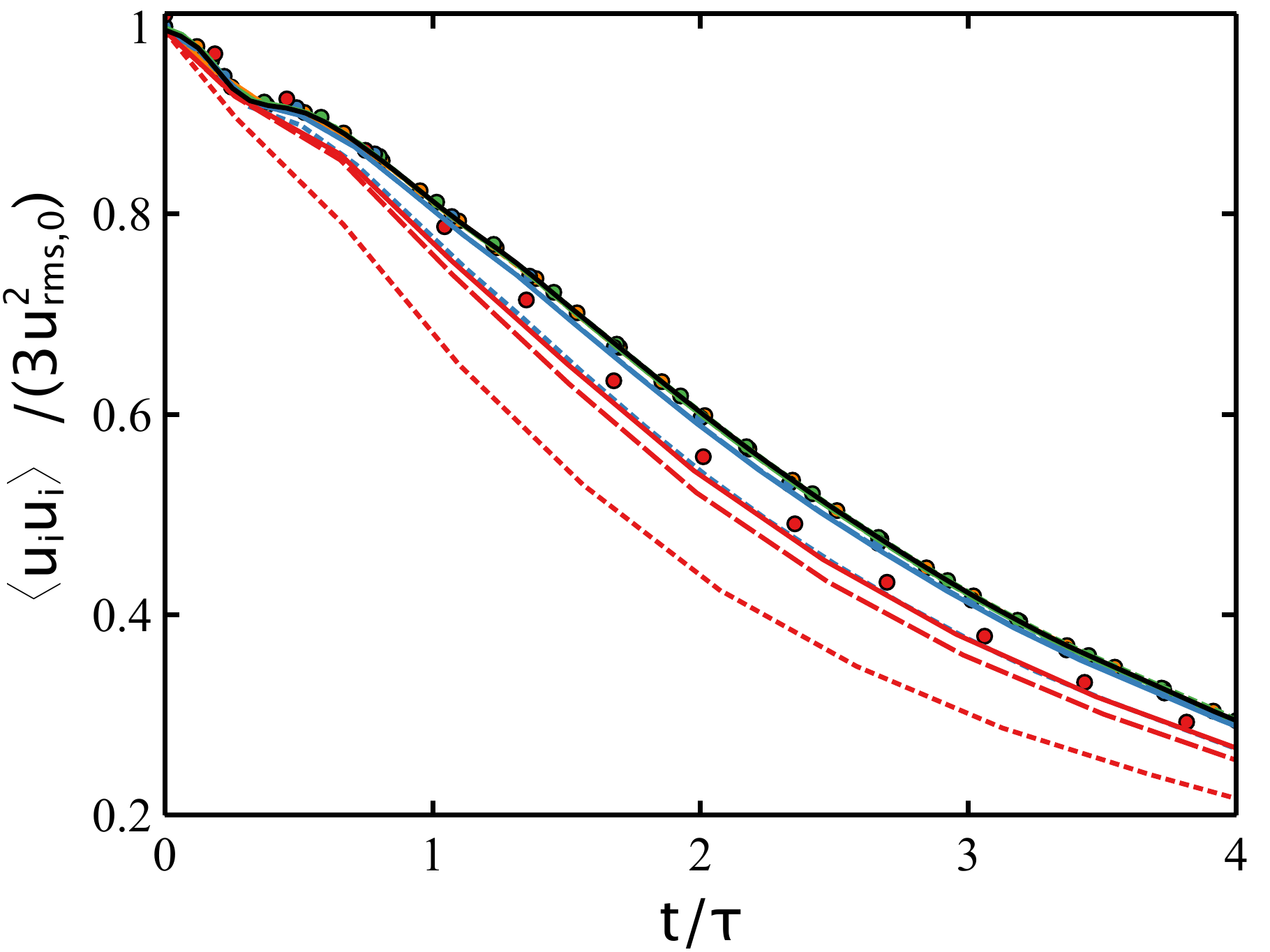}
\label{fig:HIT_temporal_a}
}
\subfloat[Enstrophy]
{
\includegraphics[width=.45\textwidth]{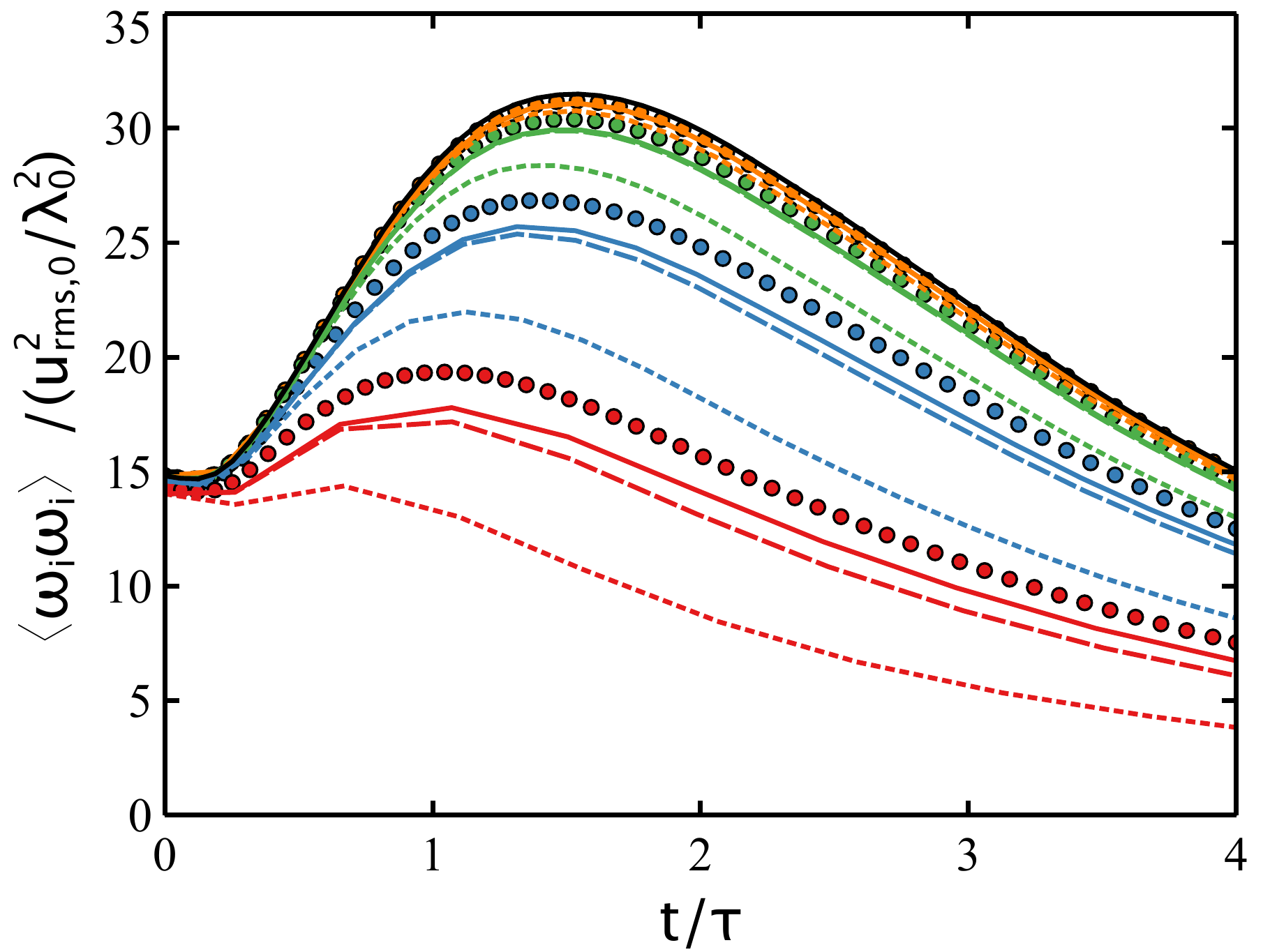}
\label{fig:HIT_temporal_b}

}
\caption{Temporal evolution of the kinetic energy and enstrophy for simulations performed with different mesh resolution and WENO-Z reconstruction methods. Legend is recalled in the text and in \cref{fig:HIT_legend}.}
\label{fig:HIT_temporal}
\end{figure}

\begin{figure}[hbt!]
\centering
\includegraphics[width=.5\textwidth]{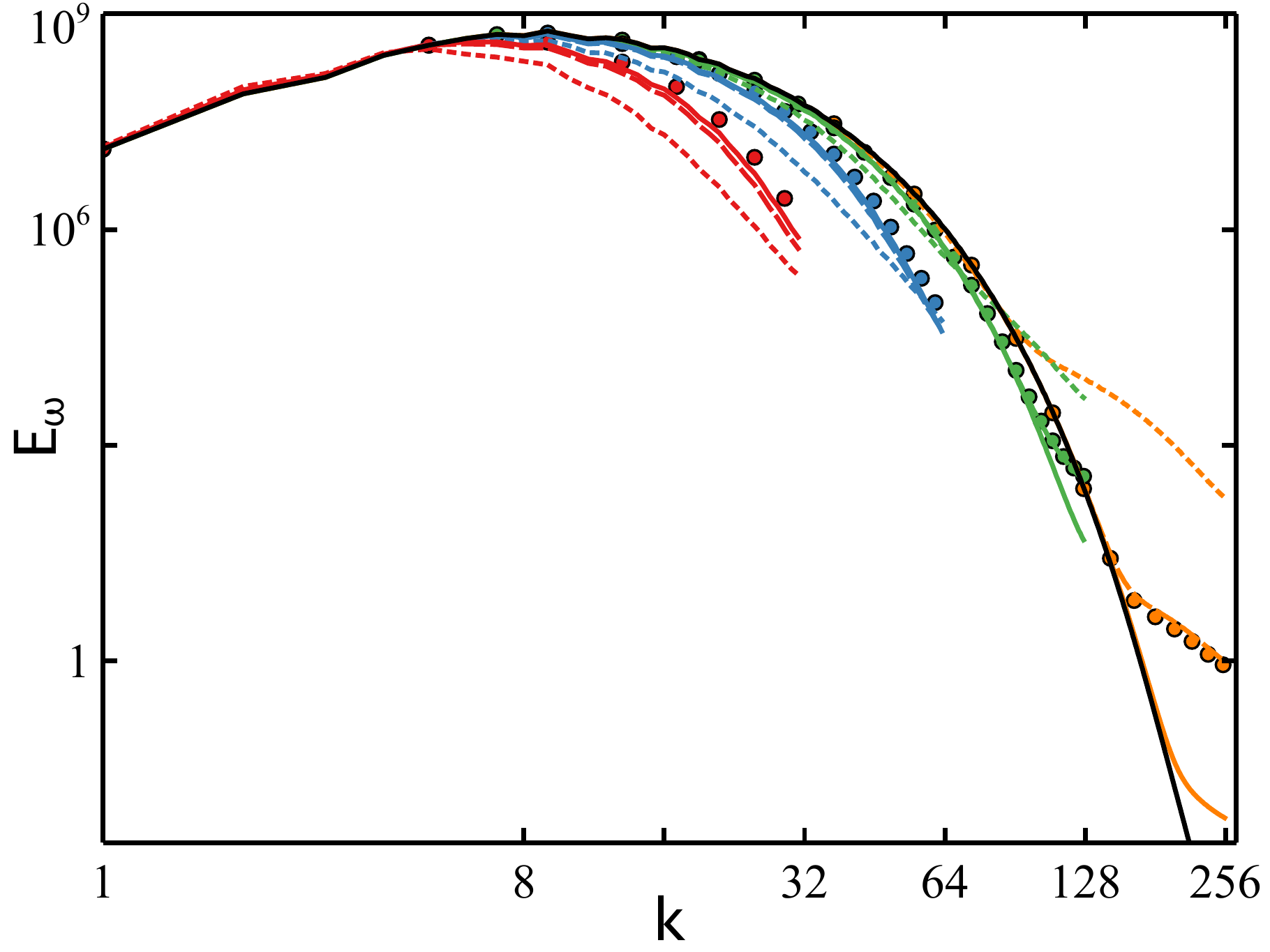}
\caption{Spectra of vorticity for simulations performed with different mesh resolution and WENO-Z reconstruction methods. Legend is recalled in the text and in \cref{fig:HIT_legend}.}
\label{fig:HIT_spectra}
\end{figure}

\begin{figure}[hbt!]
\centering
\subfloat[$N_x=64^3$]{\label{fig:HIT_spectra_WENO_zoom_a}\includegraphics[width=0.35\textwidth]{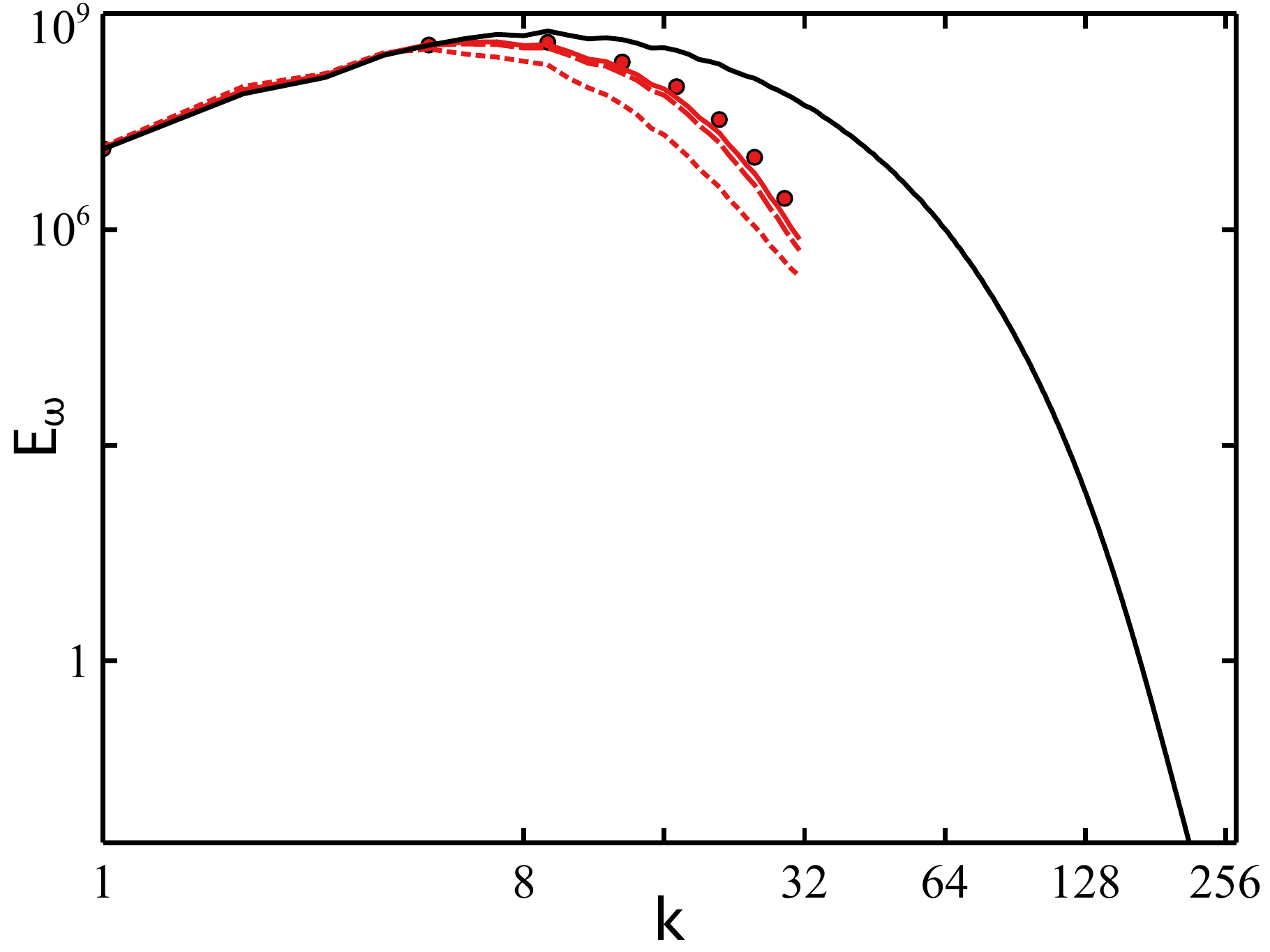}}
\subfloat[$N_x=128^3$]{\label{fig:HIT_spectra_WENO_zoom_b}\includegraphics[width=0.35\textwidth]{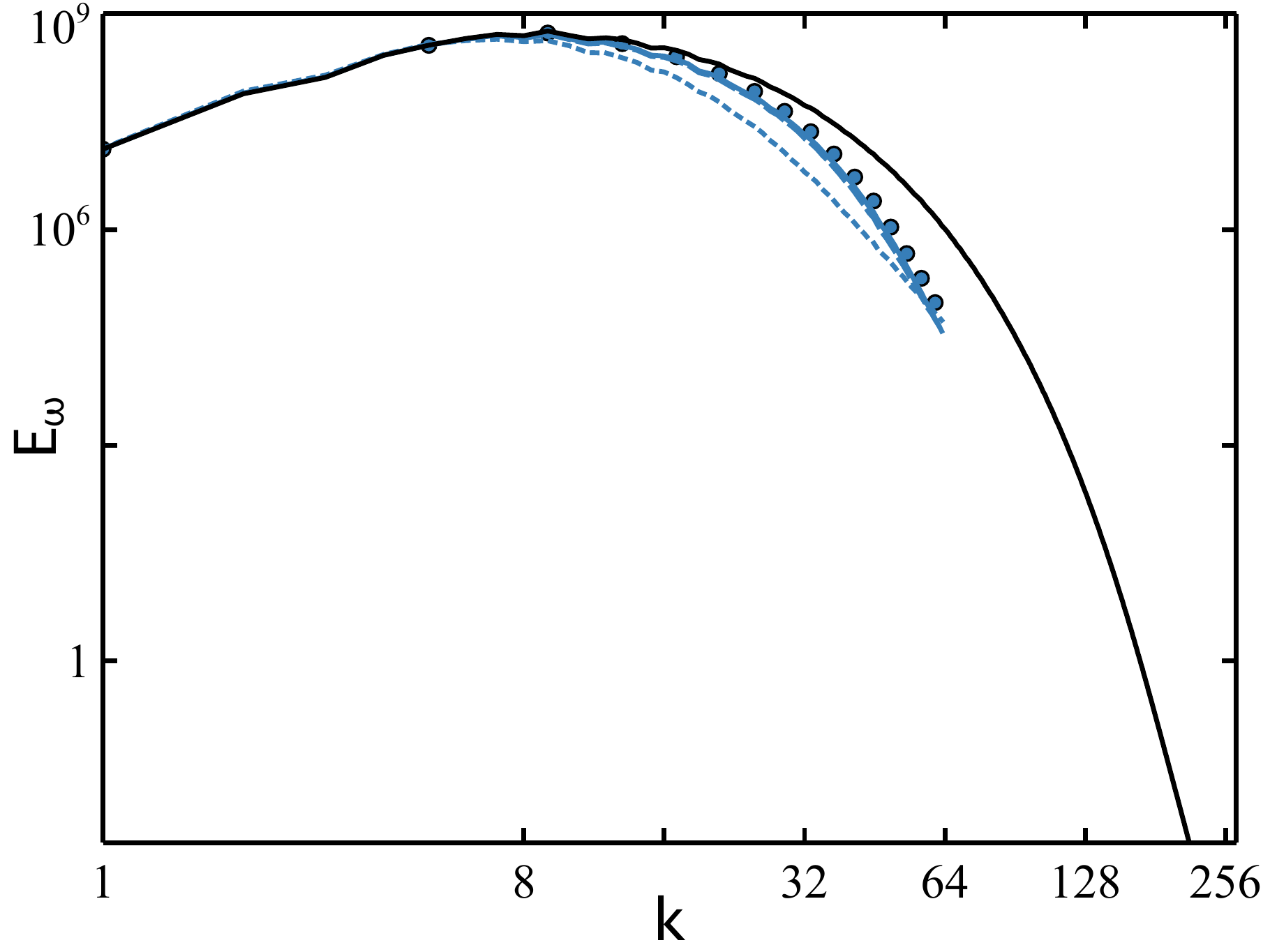}} \\
\subfloat[$N_x=256^3$]{\label{fig:HIT_spectra_WENO_zoom_c}\includegraphics[width=0.35\textwidth]{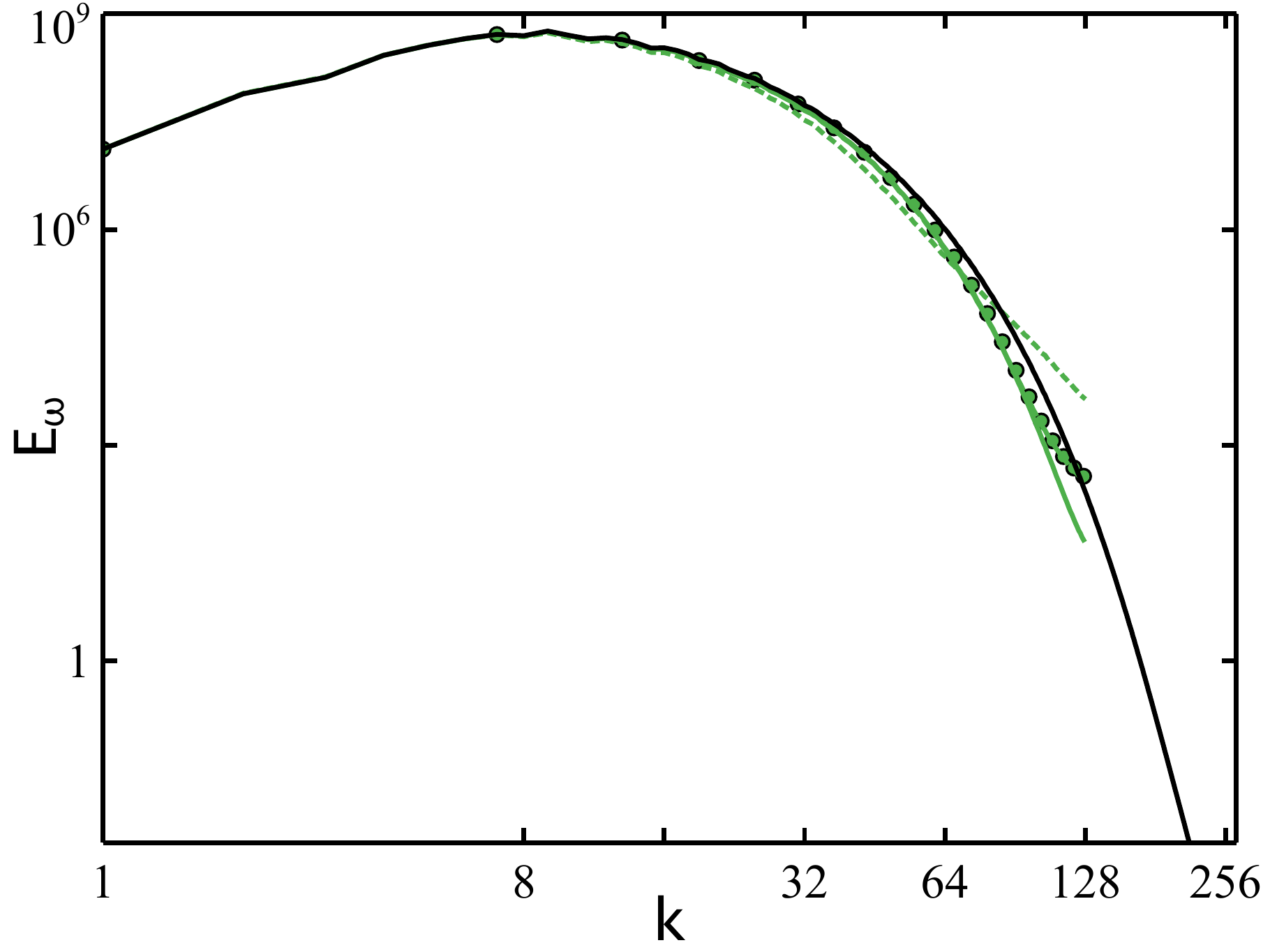}}
\subfloat[$N_x=512^3$]{\label{fig:HIT_spectra_WENO_zoom_d}\includegraphics[width=0.35\textwidth]{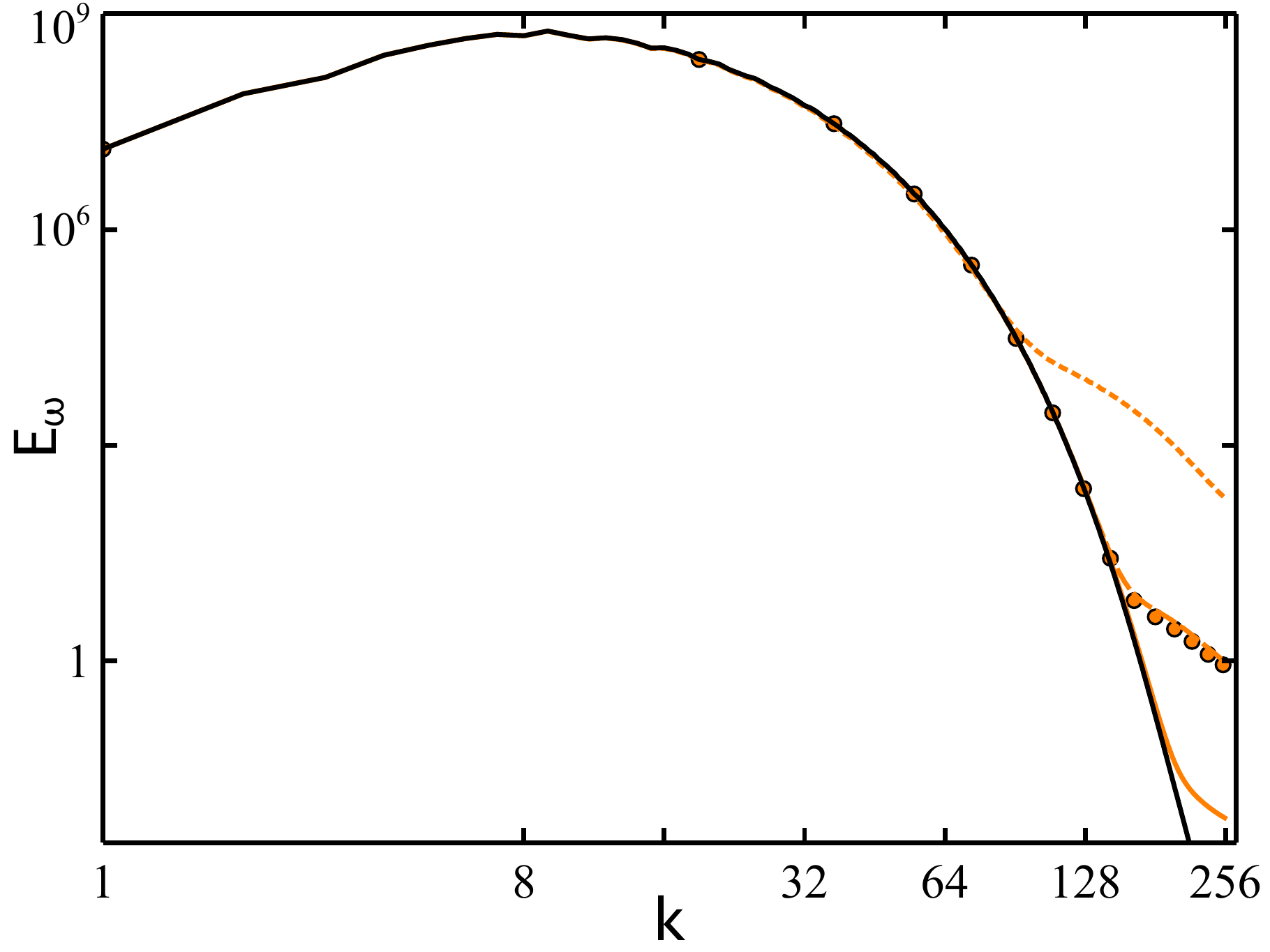}}
\caption{Extracted panels of the spectra of vorticity in \cref{fig:HIT_spectra} for results computed with different mesh resolutions.}
\label{fig:HIT_spectra_WENO_zoom}
\end{figure}

\begin{figure}[hbt!]
\centering
\includegraphics[width=.35\textwidth]{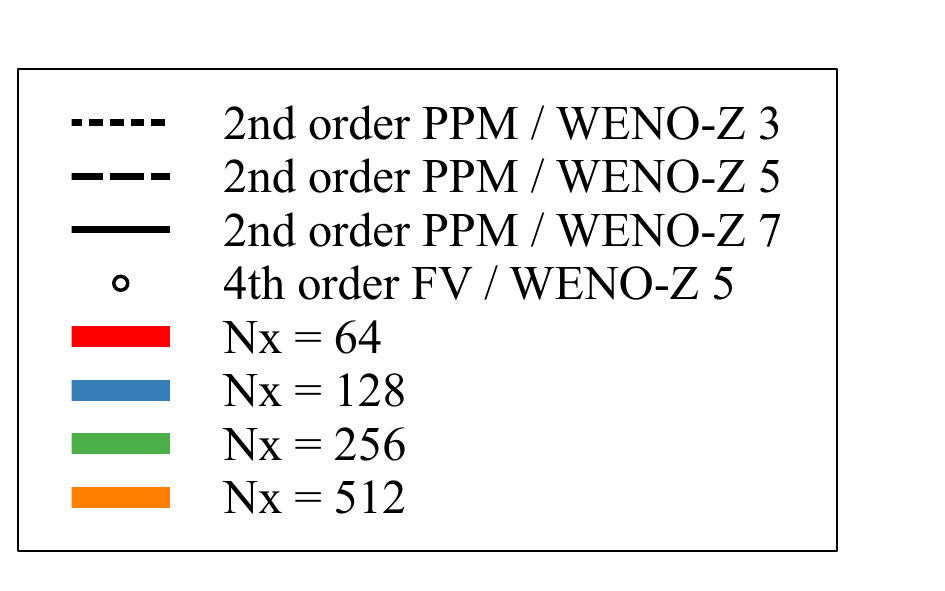}
\caption{Legend.}
\label{fig:HIT_legend}
\end{figure}

\Cref{fig:HIT_temporal} presents the temporal evolution of the kinetic energy and the enstrophy from $t=0$ to $t/\tau=4$. It can be seen that significant compressibility effects are generated quickly after the beginning of the simulation, suggesting the generation of eddy shocklets in the domain until  $t/\tau\approx 0.5$. After  $t/\tau\approx 1$, compressible shocks are no longer generated and they start to decay in a monotone manner. \Cref{fig:HIT_spectra} presents the spectra taken at $t/\tau=4$ for the vorticity. In these figures, the circle symbol represents the solution computed with the fourth-order finite-volume WENO method (\textbf{RNS} code), while the dotted, dashed and solid colored lines represent the solutions computed with the second-order hybrid PPM/WENO method (\textbf{PeleC} code) with faces data reconstructed with the third-, fifth- and seventh-order WENO-Z schemes, respectively. The red, blue, green and orange colors represent simulations performed with  $N_x=64$, $N_x=128$, $N_x=256$ and $N_x=512$, respectively. These figures contain a significant number of curves. For clarity, the legend is recalled in \cref{fig:HIT_legend} and the spectra for each mesh resolution are extracted in separated panels in \cref{fig:HIT_spectra_WENO_zoom}
for each mesh resolution. Note that the behavior of the numerical methods highlighted in \cref{fig:HIT_spectra_WENO_zoom} is virtually the same for the spectra of other physical quantities, not shown in the present paper but available in \cite{Motheau:2020} for a more detailed analysis of the hybrid PPM/WENO strategy. 

One of the most natural physical quantity to investigate when studying turbulence is the evolution of the kinetic energy. As depicted in \cref{fig:HIT_temporal_a}, for a coarse resolution ($N_x=64$) all the numerical methods over dissipate the kinetic energy over time, but it can be see that for the second-order Godunov method, the third-order WENO-Z reconstruction scheme for faces data introduces a severe dissipation. It is also interesting to note that for such coarse discretization, the fifth- and seventh-order WENO-Z schemes are able to produce results close to the fourth-order finite-volume method, whereas the whole integration scheme is theoretically second-order. When the mesh resolution is refined, virtually all the numerical methods collapse to the reference solution, at the exception of the second-order PPM method with third-order reconstruction at faces, which still provide significant extra diffusion for $N_x=128$.

It would be a mistake to only consider the kinetic energy as the main metric and to draw conclusions about the convergence of the solution respect to the discretization size. Indeed, \cref{fig:HIT_temporal_b} reveals that capturing the correct evolution of the enstrophy is challenging and that an acceptable solution is not reached below a mesh resolution of $N_x=256$. The qualitative observation of \cref{fig:HIT_temporal_a} and \cref{fig:HIT_temporal_b} can be misleading about the actual rate of convergence of the solution. In order to provide a quantitative analysis, the relative error between solutions computed with the WENO-Z schemes and the $8$th-order reference solution (\textbf{SMC} code) are evaluated for the enstrophy and the kinetic energy at $t/\tau=4$, denoted $\epsilon_{\left<\vec{\omega}_i \vec{\omega}_i \right>}$ and $\epsilon_{\left<\vec{u}_i \vec{u}_i \right>}$, respectively. The convergence rates of $\epsilon_{\left<\vec{\omega}_i \vec{\omega}_i \right>}$ and $\epsilon_{\left<\vec{u}_i \vec{u}_i \right>}$ are gathered in \cref{tab:HIT_convergence_rate_enstrophy_energy}. It can be seen that all the numerical methods, whatever their theoretical asymptotic order of accuracy, exhibit roughly a second-order rate of convergence. This is also confirmed by computing the convergence rate of the $\mathcal{L}^1$-norm of the error on the velocity at $t/\tau=4$, denoted $\varepsilon_{u}$, with results gathered in \cref{tab:HIT_convergence_rate_velocity}. The fourth-order finite-volume WENO method, as well as the second-order hybrid PPM/WENO method with reconstruction at faces with fifth- and seventh-order WENO-Z schemes exhibit an overall second-order rate of convergence. These results are consistent with the analysis performed previously at \ref{subsec:COVO} and \ref{subsec:Shu_Osher}. Indeed, in the present compressible homogeneous isotropic turbulence test case the theoretical convergence rate, valid for smooth solutions, is not recovered because of the presence of shocklets that introduce first-order errors in the solution. The more surprising result is the poor performance of the reconstruction at faces with a third-order WENO-Z scheme.

The behavior of the numerical methods can also be observed by analyzing the spectra of the turbulence. It can be seen from \cref{fig:HIT_spectra} and \cref{fig:HIT_spectra_WENO_zoom} that for  coarse mesh resolutions ($N_x = 64$ and $N_x=128$), all the numerical methods provide roughly the same results, despite the third-order WENO-Z that brings more dissipation over the spectra. However, as long as the mesh resolution is refined, each numerical method starts to exhibit a different behavior in the high frequency range. As shown in \cref{fig:HIT_spectra_WENO_zoom_c} for $N_x = 256$,  reconstruction at faces with third- and fifth-order WENO schemes present an energy pile-up in the high-frequencies. The seventh-order WENO scheme presents a better slope, closer to the reference solution. As shown in \cref{fig:HIT_spectra_WENO_zoom_d}, such behavior is exacerbated for a fine mesh resolution ($N_x = 512$). Despite the finer mesh resolution, the third-order WENO-Z scheme still starts to diverge at the wave number $k \approx 64$, and because significant errors are introduced by this scheme in the interfacial reconstruction, refinement of the mesh does not help to capture the spectra of the turbulence. The seventh-order WENO-Z scheme performs the best, with a spectra very close to the reference solution. The fifth-order WENO-Z scheme presents an intermediate solution, and interestingly enough, it can be seen that the volume integration scheme does not impact on the spectra. In other words, performing the volume integration over cells with a second-order Godunov method or a fourth-order finite-volume WENO method provide the same spectra. This suggests that in order to correctly capture the turbulent spectra, the crucial piece is the numerical method employed for the reconstruction at cell faces, not the quadrature method used for volume integration.

\begin{center}
\begin{table*}[ht]%
\caption{HIT test case: convergence rates of the relative error for the enstrophy and kinetic energy at $t/\tau=4$.}\label{tab:HIT_convergence_rate_enstrophy_energy}
\centering
\begin{tabular*}{400pt}{@{\extracolsep\fill}lccc@{\extracolsep\fill}}
\toprule
&  & \multicolumn{2}{@{}c@{}}{\textbf{Estimated Order}} \\\cmidrule{3-4}
\textbf{Software} & \textbf{Method}    & \multicolumn{1}{@{}l@{}}{\textbf{Enstropy $\left<\mathbf{\omega}_i \mathbf{\omega}_i \right>$}}  & \textbf{Kinetic Energy $\left<\mathbf{u}_i \mathbf{u}_i \right>$}   \\
\midrule
\textbf{RNS} & $4$th-order / WENO-Z 5   &  $2.26$ & $2.7$ \\
\textbf{PeleC} & Hybrid PPM/WENO-Z 3    & $2.24$  & $1.8$ \\ 
\textbf{PeleC} & Hybrid PPM/WENO-Z 5    & $2.29$  & $2.27$ \\ 
\textbf{PeleC} & Hybrid PPM/WENO-Z 7    & $2.21$  & $2.35$ \\ 
\bottomrule
\end{tabular*}
\end{table*}
\end{center}


\begin{center}
\begin{table}[hb]%
\centering
\caption{HIT test case: convergence rate of the $\mathcal{L}^1$-norm of the error on the velocity at $t/\tau=4$.}\label{tab:HIT_convergence_rate_velocity}%
\begin{tabular*}{240pt}{@{\extracolsep\fill}lcc@{\extracolsep\fill}}
\toprule
\textbf{Software} & \textbf{Method} & \textbf{Estimated Order}  \\
\midrule
\textbf{RNS} & $4$th-order / WENO-Z 5 & $2.22$    \\
\textbf{PeleC} & Hybrid PPM/WENO-Z 3 &  $1.43$    \\ 
\textbf{PeleC} & Hybrid PPM/WENO-Z 5 &  $2.05$   \\ 
\textbf{PeleC} & Hybrid PPM/WENO-Z 7 &  $1.98$    \\ 
\bottomrule
\end{tabular*}
\end{table}
\end{center}


Another aspect to take into account is the computational burden. As explained in \cite{Motheau:2020}, the major advantage of the hybrid PPM/WENO strategy is that a given similar solution can be achieved $200$ times faster than if the costly fourth-order finite-volume WENO method \cite{Titarev:2004} is employed. The main reason is that the hybrid PPM/WENO strategy only requires one WENO evaluation by face, rather than the fourth-order finite-volume that require a significant number of WENO evaluations for each point of the Gaussian quadrature rule. Indeed, as explained in \cite{Motheau:2020}, it is evaluated that for three dimensions and for only one component, achieving fourth-order accuracy requires $14$ times more interpolation procedures by cell than with the PPM algorithm. Some methods \cite{Dumbser:2007,Dumbser:2008,Boscheri:2021} have been developed to avoid integration of flux over the cell face while maintaining high order accuracy by implicit solutions for a polynomial interpolation within each cell.  However, some of the methods mentioned are partially implicit.  In \cite{Dumbser:2008} it is acknowledged that, for nonlinear fluxes, a great deal of time is spent solving these (local) problems.  In \cite{Boscheri:2021} a global implicit solve is required.  All of these methods would be interesting to implement and compare in a future study, however it is out of scope of the present paper.

In Motheau and Wakefield \cite{Motheau:2020}, only the fifth-order WENO-Z scheme was employed throughout the whole study. Here, the computational cost for each variation of accuracy of the WENO-Z scheme is studied in the context of the $N_x=64$ case. Results are reported in \cref{tab:HIT_cpu_time}. The computational cost of one WENO evaluation and the total computational time are reported in the third and last columns respectively. The $3$rd- and $5$th-order WENO-Z scheme have a very similar cost, whereas the $7$th-order scheme is slightly more costly. The second-order hybrid PPM/WENO strategy is computationally inexpensive, making use of a $7$th-order WENO-Z scheme attractive in the context of adaptive mesh refinement.

\begin{center}
\begin{table}[tbhp]%
\centering
\caption{HIT computational time for $N_x=64$}\label{tab:HIT_cpu_time}%
\begin{tabular*}{300pt}{@{\extracolsep\fill}lcc@{\extracolsep\fill}}
\toprule
\textbf{WENO-Z scheme} & \textbf{WENO evaluation time [s]} & \textbf{Total CPU time [s]}  \\
\midrule
$3$rd-order &  $0.078854$ & $134.08$  \\ 
$5$th-order &  $0.087210$ & $137.66$ \\ 
$7$th-order &  $0.10913$  & $164.41$ \\ 
\bottomrule
\end{tabular*}
\end{table}
\end{center}


\section{Conclusion}

The present paper investigates the accuracy of numerical schemes in finite-volume methods. More specifically, the impact of the third-, fifth- and seventh-order WENO-Z schemes for the reconstruction of data at cell faces has been investigated to assess their interaction with the quadrature rules employed to compute the flux over a cell face.

The accuracy of numerical methods is sometimes evaluated on canonical problems that are not representative of practical flows encountered in engineering applications. Other times, the theoretical order of convergence is viewed as the sole metric of accuracy, but the dissipation, dispersion, and computation time for finite meshes are often neglected.

Three test cases of increasing complexity have been investigated in the present work. It has been shown that on a problem with a smooth solution, the theoretical order of convergence for each method is retrieved and changing the accuracy of the scheme for the reconstruction of data at cell faces does not significantly impact the results. Furthermore, it has been demonstrated that for a shock-driven problem, all the methods collapse to first-order.

The temporal and spectral study of the decay of compressible homogeneous isotropic turbulence has revealed that using a high-order quadrature rule to compute the average over a finite-volume cell does not improve the spectral accuracy and that all methods present an actual second-order convergence for the meshes required to resolve physical scales. However it has been found that the choice of the numerical scheme to reconstruct data at cell faces is critical to correctly capture turbulent spectra.

The outcome of the present paper is that in order to achieve an accurate solution of a practical turbulent flow with finite-volume methods, it is more efficient to choose a fine mesh discretization and use a combination of a high-order reconstruction at cell faces and a low-order quadrature rule than to use a traditional high-order finite volume method.


\section*{Funding Sources}

The work here was supported by the U.S. Department of Energy, Office of Science, Office of Advanced Scientific Computing Research, Applied Mathematics program under contract number DE-AC02005CH11231. John Wakefield has been funded by The National Science Foundation (NSF) Division of Mathematical Sciences (DMS) Mathematical Sciences Graduate Internship Program, administered by the Oak Ridge Institute for Science and Education (ORISE) through an interagency agreement between the U.S. Department of Energy (DOE) and NSF. ORISE is managed by ORAU under DOE contract number DE-SC0014664. All opinions expressed in this paper are the author’s and do not necessarily reflect the policies and views of NSF, DOE or ORAU/ORISE.

\section*{Acknowledgments}
The authors would like to thank Dr. John B. Bell (Senior Scientist, Chief Scientist, Lawrence Berkeley National Laboratory) for the fruitful discussions about shock capturing numerical methods and turbulent flows.

\bibliography{bib_manu_2018}%

\end{document}